\documentclass[sn-mathphys]{sn-jnl}

\usepackage{amsfonts}
\usepackage{amsmath}
\usepackage{amssymb}
\usepackage{amsthm}
\usepackage{color}
\usepackage{commath}
\usepackage{enumerate}
\usepackage{graphics}
\usepackage{graphicx}
\usepackage{mathtools}
\usepackage[utf8x]{inputenc}

\DeclareMathOperator{\diag}{diag}

\DeclareMathOperator{\SL}{SL}

\DeclareSymbolFont{eulargesymbols}{U}{zeuex}{m}{n}
\DeclareMathSymbol{\intop}{\mathop}{eulargesymbols}{"52}
\DeclareMathSymbol{\ointop}{\mathop}{eulargesymbols}{"49}

\newcommand{\al}{\alpha}
\newcommand{\bC}{\mathbb C}

\newcommand{\bN}{\mathbb N}

\newcommand{\bR}{\mathbb R}

\newcommand{\bT}{\mathbb T}
\newcommand{\bZ}{\mathbb Z}

\newcommand{\eps}{\varepsilon}

\newcommand{\ka}{\kappa}
\newcommand{\la}{\lambda}

\newcommand{\nf}{\infty}

\newcommand{\si}{\sigma}
\newcommand{\Si}{\Sigma}

\newcommand{\tht}{\theta}

\renewcommand{\arraystretch}{1.2}

\renewcommand{\e}{\mathrm{e}}
\renewcommand{\ge}{\geqslant}
\renewcommand{\d}{\dif}

\renewcommand{\i}{\mathrm{i}}

\renewcommand{\le}{\leqslant}

\renewcommand{\x}{\raisebox{0.2mm}{\mbox{\tiny $\times$}\hspace{-0.3mm}}}

\theoremstyle{thmstyleone}
\newtheorem{theorem}{Theorem}[section]
\newtheorem{conjecture}{Conjecture}[section]

\theoremstyle{thmstyletwo}
\newtheorem{example}{Example}[section]
\newtheorem{remark}{Remark}[section]

\theoremstyle{thmstylethree}
\newtheorem{definition}{Definition}[section]
\jyear{2021}

\begin{document}

\title[Fast Toeplitz eigenvalue computations]{\Large{Fast Toeplitz eigenvalue computations, joining interpolation-extrapolation matrix-less algorithms and simple-loop theory}}

\author*[1]{\fnm{Manuel} \sur{Bogoya}}\email{johanmanuel.bogoya@uninsubria.it}

\author[2]{\fnm{Sven-Erik} \sur{Ekstr\"om}}\email{sven-erik.ekstrom@it.uu.se}

\author[1,2]{\fnm{Stefano} \sur{Serra--Capizzano}}\email{s.serracapizzano@uninsubria.it}

\affil*[1]{\orgdiv{Dipartimento di Scienza e Alta Tecnologia}, \orgname{University of Insubria}, \orgaddress{\street{via Valleggio, 11}, \city{Como}, \postcode{22100}, \country{Italy}}}

\affil[2]{\orgdiv{Department of Information Technology, Division of Scientific Computing}, \orgname{Uppsala University}, \orgaddress{\street{L\"agerhyddsv. 2}, \city{Uppsala}, \postcode{SE-751 05}, \country{Sweden}}}

\abstract{Under appropriate technical assumptions, the simple-loop theory allows to deduce various types of asymptotic expansions for the eigenvalues of Toeplitz matrices generated by a function $f$. Independently and under the milder hypothesis that $f$ is even and monotonic over $[0,\pi]$, matrix-less algorithms have been developed for the fast eigenvalue computation of large Toeplitz matrices, within a linear complexity in the matrix order: behind the high efficiency of such algorithms there are the expansions predicted by the simple-loop theory, combined with the extrapolation idea.

Here we focus our attention on a change of variable, followed by the asymptotic expansion of the new variable, and we adapt the matrix-less algorithm to the considered new setting.

Numerical experiments show a higher precision (till machine precision) and the same linear computation cost, when compared with the matrix-less procedures already presented in the relevant literature. Among the advantages, we concisely mention the following: a) when the coefficients of the simple-loop function are analytically known, the algorithm computes them perfectly; b) while the proposed algorithm is better or at worst comparable to the previous ones for computing the inner eigenvalues, it is extremely better for the computation of the extreme eigenvalues.}


\maketitle

\section{Introduction}

The target of this note is to design fast procedures for the computation of all the spectra of large Toeplitz matrices having an even generating function which is monotone in the interval $(-\pi,\pi)$. For the formal definition of Toeplitz matrix generated by a Lebesgue integrable function over the basic interval $Q\equiv(-\pi,\pi)$ see the first lines of Section \ref{sc:prel}.

This topic has been studied in the recent years by several researchers. Indeed, taking into account a clear numerical evidence developed in a systematic series of the numerical tests, in \cite{EkGa18} the second author formulated the following conjecture.
\begin{conjecture}\label{conj}
Let $l,g$ be two real-valued even functions with $g>0$ on $(0,\pi)$, and suppose that $f\equiv \frac{l}{g}$ is monotone increasing over $(0,\pi)$. Set $X_{n}\equiv T_{n}^{-1}(g)T_{n}(l)$ for all $n$. Then, for every integer $\al\ge0$, every $n$, and every $j=1,\ldots,n$, the following asymptotic expansion holds:
\begin{equation*}
\la_{j}(X_{n})=f(\tht_{j,n})+\sum_{k=1}^{\al-1} c_{k}(\tht_{j,n})h^{k}+E_{j,n,\al},
\end{equation*}
where:
\renewcommand\labelitemi{\tiny$\bullet$}
\begin{itemize}
\item the eigenvalues of $X_{n}$ are arranged in non-decreasing order, that is\\$\la_{1}(X_{n})\le\cdots\le\la_{n}(X_{n})$;\,\footnote{Note that the eigenvalues of $X_{n}$ are real, because $T_{n}(g)$ is symmetric positive definite and $X_{n}$ is similar to the symmetric matrix $T_{n}^{-\frac{1}{2}}(g)T_{n}(l)T_{n}^{-\frac{1}{2}}(g)$.}
\item $\{c_{k}\}_{k=1}^{\nf}$ is a sequence of functions from $(0,\pi)$ to $\bR$ which depends only on $l,g$;
\item $h\equiv\frac{1}{n+1}$ and $\tht_{j,n}\equiv\frac{j\pi}{n+1}=j\pi h$;
\item $E_{j,n,\al}=O(h^{\al})$ is the remainder (the error), which satisfies the inequality $\abs{E_{j,n,\al}} \le \ka_{\al} h^{\al}$ for some constant $\ka_{\al}$ depending only on $\al,l,g$.
\end{itemize}
\end{conjecture}
\noindent In the case where $g=1$ identically, Conjecture~\ref{conj} was originally formulated and supported through numerical experiments in \cite{EkGa18}. Then the algorithmic proposal was extended and refined in \cite{AhAl18,EkFu18a,EkFu18b,EkGa19}. 
When $g\equiv 1$ and $l$ satisfies further technical additional assumptions, those of the simple-loop method, Conjecture~\ref{conj} was formally proved by Bogoya, B\"ottcher, Grudsky, and Maximenko in a series of papers \cite{BoBo15a,BoBo16,BoBo17}.

However, again in the case of $g\equiv1$, the power of the simple-loop method has been not exploited completely, since in the case of a continuous generating function, the distribution results reported in Theorem \ref{teoszego-tyr} imply that $\la_{j}(T_{n}(f))=f(s_{j,n})$ with $s_{j,n}$ belonging to $(0,\pi)$ and distributed as the identity function. However, more is known and indeed also in the independent variable $\tht$ there exists an asymptotic expansion regarding exactly the points $s_{j,n}$.

This note deals with the adaptation of the interpolation-extrapolation algorithms to the previous change of variable, joined with a trick at the end points introduced in \cite{BoSe21}. 

The numerical results are extremely precise, even compared with the already good performances described in \cite{AhAl18,EkFu18a,EkFu18b,EkGa19,EkGa18}, since it is not difficult to reach machine precision, and the complexity is still linear.

The present note is organized as follows. Preliminary definitions, tools, and results are concisely reported in Section \ref{sc:prel}. Section~\ref{sc:algo} presents the new adapted algorithm for computing the Toeplitz eigenvalues: as in \cite{EkGa19}, our technique combines the extrapolation procedure proposed in \cite{AhAl18,EkGa18} -- which allows the computation of {\em some} of the eigenvalues of $X_{n}$ -- with an appropriate interpolation process, designed for the simultaneous computation of {\em all} the eigenvalues of $X_{n}$, with the additional end point trick in \cite{BoSe21}. In Section~\ref{sc:num} we present the numerical experiments, while in Section \ref{conclusions} we draw conclusions and we list few open problems for research lines to be investigated in the next future. 

\section{Preliminaries and Tools}\label{sc:prel}

For a real or complex valued function $f$ in $L^{1}[-\pi,\pi]$, let $\mathfrak{a}_{j}(f)$ be its $j$th Fourier coefficient, i.e.
\[
\mathfrak{a}_{j}(f)\equiv\frac{1}{2\pi}\int_{-\pi}^{\pi} f(\tht)\e^{-\i j\tht}\dif\tht,\quad j\in\bZ,
\]
and consider the sequence $\{T_{n}(f)\}_{n=1}^{\nf}$ of the $n\times n$ Toeplitz matrices defined by $T_{n}(f)\equiv(\mathfrak{a}_{j-k}(f))_{j,k=0}^{n-1}$. The function $f$ is customarily referred to as the generating function of this sequence.

As a second step, we introduce some notations and definitions concerning general sequences of matrices. For any function $F$ defined on the complex field and
for any matrix $A_{n}$ of size $d_{n}$, by the symbol $\Si_{\la}(F,A_{n})$, we denote the mean
\begin{equation*}
\Si_{\la}(F,A_{n})=\frac{1}{d_{n}} \sum_{j=1}^{d_{n}}
F[\la_{j}(A_{n})],
\end{equation*}
 while by the symbol $\Si_{\si}(F,A_{n})$, we denote the mean
\begin{equation*}
\Si_{\si}(F,A_{n})=\frac{1}{d_{n}} \sum_{j=1}^{d_{n}}
F[\si_{j}(A_{n})].
\end{equation*}
\begin{definition}
Given a sequence $\{A_{n}\}$ of matrices of size $d_{n}$ with $d_{n}<d_{n+1}$ and given a Lebesgue-measurable function $\psi$ defined over a measurable set $K\subset {\bR}^{\nu}$, and $\nu \in \bN^{+}$, of finite and positive Lebesgue measure $\mu(K)$, we say that $\{A_{n}\}$ is distributed as $(\psi,K)$ in the sense of the eigenvalues if for any continuous function $F$ with bounded support, the following limit relation holds
\begin{equation*}
\lim_{n\rightarrow \nf}\Si_{\la}(F,A_{n})= \frac{1}{\mu(K)}\int_{K} F(\psi)\,\d\mu.
\end{equation*}
In this case, we write in short $\{A_{n}\}\sim_{\la} (\psi,K)$. Furthermore we say that $\{A_{n}\}$ is distributed as $(\psi,K)$ in the sense of the singular values if for any continuous function $F$ with bounded support, the following limit relation holds
\begin{equation*}
\lim_{n\rightarrow \nf}\Si_{\si}(F,A_{n})= \frac{1}{\mu(K)}\int_{K} F(\,\abs{\psi})\,\d\mu.
\end{equation*}
In this case, we write in short $\{A_{n}\}\sim_{\si} (\psi,K)$, which is equivalent to say that $\{A_{n}^*A_{n}\}\sim_{\la} (\,\abs{\psi}^2,K)$.
\end{definition}

In Remark~\ref{rem:meaning-distribution} we provide an informal meaning of the notion of eigenvalue distribution. For the singular value distribution similar statements can be written.

\begin{remark}\label{rem:meaning-distribution}
The informal meaning behind the above definition is the following. If $\psi$ is continuous, $n$ is large enough, and
\begin{equation*}
\Big\{{\bf x}_{j}^{(d_{n})},\ j=1,\ldots, d_{n}\Big\}
\end{equation*}
is an equispaced grid on $K$, then a suitable ordering $\la_{j}(A_{n})$, $j=1,\ldots,d_{n}$, of the eigenvalues of $A_{n}$ is such that the pairs $\big\{\big({\bf x}_{j}^{(d_{n})},\la_{j}(A_{n})\big),\ j=1,\ldots,d_{n}\big\}$ reconstruct approximately the hypersurface
\begin{equation*}
\big\{({\bf x},\psi({\bf x})),\ {\bf x}\in K\big\}.
\end{equation*}
 In other words, the spectrum of $A_{n}$ `behaves' like a uniform sampling of $\psi$ over $K$.
For instance, if
$\nu=1$, $d_{n}=n$, and $K=[a,b]$, then the eigenvalues of $A_{n}$ are approximately equal to $\psi\big(a+\frac{j}{n}(b-a)\big)$, $j=1,\ldots,n$, for $n$ large enough and up to at most $o(n)$ outliers.
Analogously, if $\nu=2$, $d_{n}=n^2$, and $K=[a_{1},b_{1}]\times [a_2,b_2]$, then the eigenvalues of $A_{n}$ are approximately equal to $\psi\big(a_{1}+\frac{j}{n}(b_{1}-a_{1}),a_2+\frac{k}{n}(b_2-a_2)\big)$, $j,k=1,\ldots,n$, for $n$ large enough and up to at most $o(n^2)$ outliers.
\end{remark}

The asymptotic distribution of eigenvalues and singular values of Toeplitz matrix sequences has been studied deeply and continuously in the last century (for example see \cite{BaGa20b,BaGa20a,BoSi99,GaSe17,GaSe18} and references therein).
The starting point of this theory, which contains many extensions and other results, is a famous theorem of Szeg\H{o}~\cite{GrSZ84}, which we report in 
the version due to Tyrtyshnikov and Zamarashkin~ \cite{TyZa98}.
\begin{theorem}
\label{teoszego-tyr}
If $f$ is integrable over $Q\equiv(-\pi,\pi)$, and if $\{T_{n}(f)\}$ is the sequence of
Toeplitz matrices generated by $f$, then
\begin{equation*}
\{ T_{n}(f)\}\sim_{\si} (f,Q).
\end{equation*}
Moreover, if $f$ is also real-valued, then each matrix $T_{n}(f)$
is Hermitian and
\begin{equation*}
\{T_{n}(f)\}\sim_{\la} (f,Q).
\end{equation*}
\end{theorem}
Furthermore, strong localization results are known in the case where the generating function is real-valued, as stated in \cite[Th.2.2]{Se97b} which we partly report below.
\begin{theorem}
\label{strong-loc}
If $f$ is integrable and real-valued almost everywhere over $Q\equiv(-\pi,\pi)$, $m$ is the essential infimum of $f$, $M$ is the essential supremum of $f$, and if $\{T_{n}(f)\}$ is the sequence of Toeplitz matrices generated by $f$, then
\begin{equation*}
\lambda_{j}(T_n(f)) \in (m,M)
\end{equation*}
for every $j=1,\ldots,n$, for every positive integer $n$, under the assumption $m<M$.
Moreover, if $m=M$ then the generating function $f$ is constant almost everywhere and trivially we conclude that $T_{n}(f)$ coincides with $m$ times the identity matrix. 
\end{theorem}

\subsection{The Simple-Loop case}

For $\al>0$, the well-known weighted Wiener algebra $W^{\al}$ is the collection of all functions $f\colon\bT\to\bC$ whose Fourier coefficients satisfy
\[\norm{f}_{\al}\equiv\sum_{j=-\nf}^{\nf}\abs{\mathfrak{a}_{j}(f)}(\,\abs{j} +1)^{\al}<\nf.\]
It is easy to see that if $f\in W^{\al}$ then $f\in C^{\lfloor\al\rfloor}[-\pi,\pi]$, hence the constant $\al$ is measuring the smoothness of the symbol $f$. In what follows we extend every symbol $f$ to the whole real line in the natural way turning it into a $2\pi$-periodic function, and we denote this extension by $f$ as well.

The simple-loop class, denoted by $\SL^{\al}$, consists in the collection of all the real-valued symbols in $W^{\al}$ tracing out a simple-loop over the interval $[-\pi,\pi]$ with the following properties:
\begin{enumerate}[(i)]
\item the range of $f$ is a segment $[0,\mu]$ with $\mu>0$;
\item $f(0)=f(2\pi)=0$, $f''(0)=f''(2\pi)>0$;
\item there is a unique $\tht_{0}\in[0,2\pi]$ such that $f(\tht_{0})=\mu$, $f'(\tht)>0$ for $\tht\in(0,\tht_{0})$ and $f'(\tht)<0$ for $\tht\in(\tht_{0},2\pi)$.
\end{enumerate}
Note that if $f\in\SL^{\al}$ then $f'(0)=f'(2\pi)=0$, and that $\tht_{0}=\pi$ for every even symbol, i.e. a symbol satisfying $f(\tht)=f(-\tht)$ for $\tht\in\bR$.

Consider an even symbol $f\in\SL^{\al}$ with $\al>2$. For a $n\times n$ matrix $A$ let $\la_{j}(A)$ $(j=1,\ldots,n)$ be its eigenvalues. In the works \cite{BoBo15a,BoBo16,BoGr17,BoSe21}, for example, the authors state that
\begin{enumerate}[(i)]
\item the eigenvalues of $T_{n}(f)$ are all distinct, i.e.
\[\la_{1}(T_{n}(f))<\la_{2}(T_{n}(f))<\cdots<\la_{n}(T_{n}(f));\]
\item the numbers $s_{j,n}\equiv \big[f\mid_{[0,\pi]}\big]^{-1}(\la_{j}(T_{n}(f)))$ $(j=1,\ldots,n)$ satisfy
\begin{equation}\label{eq:Magic}
(n+1)s_{j,n}+\eta(s_{j,n})=\pi j+E_{j,n,\al},
\end{equation}
where $E_{j,n,\al}$ is an error term satisfying certain bounding condition, and $\eta$ is a function depending only on $f$ with certain smoothness depending on $\al$.
\item the previous equation \eqref{eq:Magic} has exactly one solution $s_{j,n}\in[0,\pi]$ for each $j=1,\ldots,n$.
\end{enumerate}
In the bulk of the so-called \textit{simple-loop method}, the authors use the Banach fixed-point theorem to iterate over \eqref{eq:Magic} and solve it for $s_{j,n}$, obtaining an expansion of the kind (see \cite[Th.2.2]{BoBo15a} for example)
\begin{equation}\label{eq:MainExp}
\la_{j}(T_{n}(f))= f(s_{j,n}), \qquad s_{j,n}=\tht_{j,n}+\sum_{k=1}^{\lfloor\al\rfloor}r_{k}(\tht_{j,n})h^{k}+E_{j,n,\al},
\end{equation}
where
\renewcommand\labelitemi{\tiny$\bullet$}
\begin{itemize}
\item the numbers $s_{j,n}$ are arranged in nondecreasing order;
\item $h\equiv\frac{1}{n+1}$ and $\tht_{j,n}\equiv \pi j h$;
\item the coefficients $r_{k}$ depend only on $f$ and can be found explicitly, for example
\begin{center}
\begin{tabular}{ll}
$r_{1}=-\eta$, & $r_{3}=-\eta(\eta')^{2}-\frac{1}{2}\eta^{2}\eta''$,\\
$r_{2}=\eta \eta'$, & $r_{4}=\eta(\eta')^{3}+\frac{3}{2}\eta^{2}\eta'\eta''+\frac{1}{6}\eta^{3}\eta'''$.
\end{tabular}
\end{center}
\item $E_{j,n,\al}=O(h^{\al})$ is the remainder (error) term, which satisfies the bounding $\abs{E_{j,n,\al}}\le\ka_{\al} h^{\al}$ for some constant $\ka_{\al}$ depending only on $\al$ and $f$.
\end{itemize}
Using the expansion \eqref{eq:MainExp} and the smoothness of the symbol $f$, the authors apply $f$ to both sides obtaining an expansion of the kind
\begin{equation}\label{eq:LaExp}
\la_{j}(T_{n}(f))=f(\tht_{j,n})+\sum_{k=1}^{\lfloor\al\rfloor} c_{k}(\tht_{j,n})h^{k}+E_{j,n,\al},
\end{equation}
with similar characteristics but where the coefficients $c_{k}$ involve the symbol and its derivatives, for instance $c_{1}=-f'\eta$ and $c_{2}=f'\eta\eta'+\frac{1}{2}f''\eta^{2}$. The previous works \cite{BoSe21,EkFu18b,EkGa19,EkGa18} used \eqref{eq:LaExp} as the basic expansion. We noticed that \eqref{eq:LaExp} is absorbing all the almost non-increasining consequences of $f$ at $\{0,\pi\}$ and all the related troubles that the derivatives of $f$ can produce. Hence, we decided to work with \eqref{eq:MainExp} instead.

Indeed, in the light of Theorem \ref{teoszego-tyr}, since $f$ is even and real-valued, we find $\{T_{n}(f)\}\sim_{\la} (f,[0,\pi])$. Consequently we deduce
\[
\{{\diag}_{j=1,\ldots,n}(s_{j,n})\}\sim_{\la} ({\rm id},[0,\pi]),\ \ \ \ {\rm id}(\tht)\equiv \tht,
\]
with $s_{j,n}\in (0,\pi)$, by virtue of Theorem \ref{strong-loc}.
Notice that the function id$(\cdot)$ is very basic, and, as already claimed, when compared with the study and the proposals in \cite{BoBo15a,BoBo16,BoBo17,BoSe21,EkFu18a,EkFu18b,EkGa19,EkGa18}, the troubles that the derivatives of $f$ produce are completely removed: these nice features are clearly evident by looking at the high precision of the numerical computations, reported in Section \ref{sc:num}, containing the numerical experiments.

\section{The Algorithm}\label{sc:algo}

\begin{figure}[ht]
\centering
\includegraphics[width=0.97\textwidth]{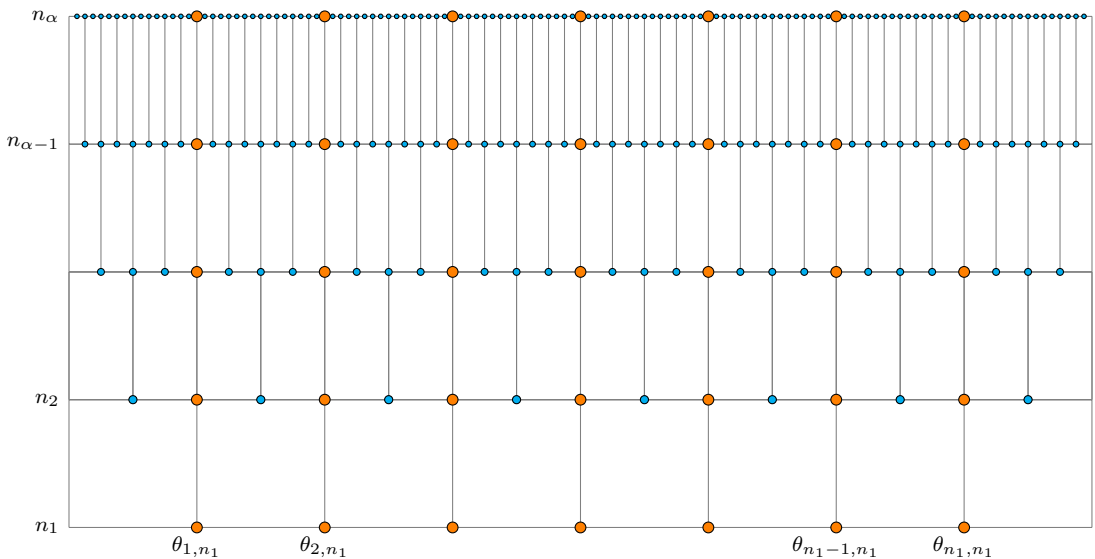}
\caption{The regular grids $\{\tht_{j,n_{k}}\}$ for $j=1,\ldots,n_{k}$ and $k=1,\ldots,\al$. In the precomputing phase we will need to calculate the eigenvalues of $T_{n_{k}}(f)$ for $k=1,\ldots,\al$ corresponding to the blue and orange dots combined, but in the interpolation phase, we will use only the eigenvalues corresponding to the orange dots, that is $\la_{j_{k}}(T_{n_{k}}(f))$ for $j_{1}=1,\ldots,n_{1}$ and $k=1,\ldots,\al$.}\label{fg:Grid}
\end{figure}

Our algorithm is based on the expansion \eqref{eq:MainExp} and is an evolution of the algorithms proposed in \cite{BoSe21,EkFu18b,EkGa19,EkGa18}. As in the mentioned works it is suited for parallel implementation and can be called matrix-less since it does not require to calculate or even, to store the matrix entries. For every $n\in\bN$ let $h\equiv\frac{1}{n+1}$ and $\tht_{j,n}\equiv \pi jh$, thus the collection $\{\tht_{j,n}\}_{j=0}^{n+1}$ is a regular grid for the interval $[0,\pi]$ with step size $\pi h$. We will use similar notations for $n$ and $h$, for example, $h_{k}$ means $\frac{1}{n_{k}+1}$, and so on. We assume that
\renewcommand\labelitemi{\tiny$\bullet$}
\begin{itemize}
\item the symbol $f$ is even and real-valued, strictly increasing in the interval $[0,\pi]$, and $f(0)=0$;
\item $n_{1}$ and $\al$ are fixed natural numbers and $n\gg n_{1}$;
\item for $k=1,\ldots,\al$ let $n_{k}\equiv 2^{k-1}(n_{1}+1)-1$;
\item for $j_{1}=1,\ldots,n$ and $k=1,\ldots,\al$, let $j_{k}\equiv 2^{k-1}j_{1}$.
\end{itemize}
Note that $j_{k}$ depends on $j_{1}$, and similarly, $n_{k}$ depends on $n_{1}$, but for notation simplicity, we suppressed those dependencies.
The indexes $j_{k}$ and the matrix sizes $n_{k}$ were calculated in such a way that
\[\tht_{j_{1},n_{1}}=\tht_{j_{2},n_{2}}=\cdots=\tht_{j_{\al},n_{\al}},\]
which is the key idea of the following extrapolation phase, see Figure \ref{fg:Grid}.
We also want to emphasize that $f$ is not necessarily simple-loop.

As in \cite{BoSe21,EkFu18b,EkGa19,EkGa18}, our algorithm is designed to calculate eigenvalues for ``big'' matrix sizes $n$ with respect to $n_{1},\ldots,n_{\al}$, meaning that, from a computational viewpoint, the calculation of the eigenvalues of $T_{n}(f)$ is hard while for $T_{n_{k}}(f)$ can be easily done with any standard eigensolver (i.e. \texttt{Eigenvalues} in \textsc{Mathematica} or \texttt{eig} in \textsc{Matlab}). But our proposal is able to reach machine precision accuracy easily. The algorithm has two phases, the first one involves an extrapolation procedure, and the second one consists in a local interpolation technique.
As a precomputing phase we need to calculate the eigenvalues of $T_{n}(f)$ for $n=n_{1},\ldots,n_{\al}$.

\noindent{\textbf{Extrapolation} For each fixed $j_{1}=1,\ldots,n_{1}$ let $\si_{j_{1}}\equiv\tht_{j_{1},n_{1}}=\cdots=\tht_{j_{\al},n_{\al}}$ (see the orange dots in Figure~\ref{fg:Grid}), and apply $\al$ times the expansion \eqref{eq:MainExp} obtaining
{\small\begin{eqnarray*}
\la_{j_{1}}(T_{n_{1}}(f))-f(\si_{j_{1}})&=&r_{1}(\si_{j_{1}})h_{1}+r_{2}(\si_{j_{1}})h_{1}^{2}+\cdots+r_{\al}(\si_{j_{1}})h_{1}^{\al}+E_{j_{1},n_{1},\al},\\
\la_{j_{2}}(T_{n_{2}}(f))-f(\si_{j_{1}})&=&r_{1}(\si_{j_{1}})h_{2}+r_{2}(\si_{j_{1}})h_{2}^{2}+\cdots+r_{\al}(\si_{j_{1}})h_{2}^{\al}+E_{j_{2},n_{2},\al},\\
&\vdots&\\
\la_{j_{\al}}(T_{n_{\al}}(f))-f(\si_{j_{1}})&=&r_{1}(\si_{j_{1}})h_{\al}+r_{2}(\si_{j_{1}})h_{\al}^{2}+\cdots+r_{\al}(\si_{j_{1}})h_{\al}^{\al}+E_{j_{\al},n_{\al},\al}.
\end{eqnarray*}}
Let $\hat r_{k}(\si_{j_{1}})$ be the approximation of $r_{k}(\si_{j_{1}})$ obtained by removing all the error terms $E_{j_{k},n_{k},\al}$ and solving the resulting linear system:
{\small\begin{equation}\label{eq:IntPhase}
\begin{bmatrix}
h_{1} & h_{1}^{2} & \cdots & h_{1}^{\al}\\
h_{2} & h_{2}^{2} & \cdots & h_{2}^{\al}\\
\vdots & \vdots & \ddots & \vdots\\
h_{\al} & h_{\al}^{2} & \cdots & h_{\al}^{\al}
\end{bmatrix}
\begin{bmatrix}
\hat r_{1}(\si_{j_{1}}) \\
\hat r_{2}(\si_{j_{1}}) \\
\vdots\\
\hat r_{\al}(\si_{j_{1}})
\end{bmatrix}
=
\begin{bmatrix}
\la_{j_{1}}(T_{n_{1}}(f)) \\
\la_{j_{2}}(T_{n_{2}}(f)) \\
\vdots\\
\la_{j_{\al}}(T_{n_{\al}}(f))
\end{bmatrix}
-f(\si_{j_{1}})
\begin{bmatrix}
1\\ 1\\ \vdots\\ 1
\end{bmatrix}.
\end{equation}}
As mentioned in previous works, a variant of this extrapolation strategy was first suggested by Albrecht Böttcher in \cite[Sc.7]{BoBo15a} and is analogous to the Richardson extrapolation employed in the context of Romberg integration \cite[Sc.3.4]{StBu10}.

\noindent\textbf{Interpolation} For any $j\in\{1,\ldots,n\}$ we will estimate $r_{k}(\tht_{j,n})$. If $\tht_{j,n}$ coincides with one of the points in the grid $\{\tht_{0,n_{1}},\ldots,\tht_{n_{1}+1,n_{1}}\}$, then we have the approximations $\hat r(\tht_{j,n})$ from the extrapolation phase for free. In any other case, we will do it by interpolating the data
\[(\tht_{0,n_{k}},\hat r_{k}(\tht_{0,n_{k}})),(\tht_{1,n_{k}},\hat r_{k}(\tht_{1,n_{k}})),\ldots,(\tht_{n_{k}+1,n_{k}},\hat r_{k}(\tht_{n_{k}+1,n_{k}})),\]
for $k=1,\ldots,\al$, and then evaluating the resulting polynomial at $\tht_{j,n}$. This interpolation can be done in many ways, but to avoid spurious oscillations explained by the Runge phenomenon \cite[p.78]{Da75a}, and following the strategy of the previous works, we decided to do it considering only the $\al-k+5$ points in the grid $\{\tht_{0,n_{1}},\ldots,\tht_{n_{1}+1,n_{1}}\}$ which are closest to $\tht_{j,n}$. Those points can be determined uniquely unless $\tht_{j,n}$ is the mid point of two consecutive points in the grid, in which case we can take any of the possible two choices.

Finally, our eigenvalue approximation with $k$ terms, is given by
\begin{equation}\label{eq:NAS}
\hat\la_{j,k}^{\textrm{NAS}}(T_{n}(f))\equiv f\Big(\tht_{j,n}+\sum_{\ell=0}^{k-1} \hat r_{\ell}(\tht_{j,n})h^{\ell}\Big),
\end{equation}
where $k=1,\ldots,\al$, and NAS stands for ``Numerical Algorithm in the variable $s_{j,k}$''.

\begin{remark}
To get the best result possible, in the precomputing phase, we advise to use the proposed algorithm \eqref{eq:NAS}, calculating the eigenvalues of the matrices $T_{n}(f)$ for $n=n_{1},\ldots,n_{\al}$ with a significant number of precision digits, let's say $60$. For $\al=5$ and $n_{1}=100$, for instance, this can be done in a standard computer in a few minutes, and it only needs to be done once. While for the extrapolation phase \eqref{eq:IntPhase}, we advice to do it with the largest $\al$ possible, even if less terms will be used.
\end{remark}

\section{Numerical Experiments}\label{sc:num}

Let $\la_{j,k}^{\SL}(T_{n}(f))$ be the $k$th term approximation of $\la_{j}(T_{n}(f))$ obtained with the simple-loop method in the variable $s_{j,n}$, that is
\[\la_{j,k}^{\SL}(T_{n}(f))\equiv
f\Big(\tht_{j,n}+\sum_{k=1}^{k-1} r_{k}(\tht_{j,n})h^{k}\Big),\]
where $\tht_{j,n}=\pi jh$. Let also $\la_{j,k}^{\textrm{NA}}(T_{n}(f))$ be the $k$th term approximation of $\la_{j}(T_{n}(f))$ given by the Numerical Algorithm in \cite{EkGa19} and finally, let $\la_{j,k}^{\textrm{MNA}}(T_{n}(f))$ be the respective approximation given by the Modified Numerical Algorithm \cite[Sc.4]{BoSe21}. In order to compare the results of the different methods we use the following notation for the absolute individual errors

{\small\begin{align*}
\eps^{\SL}_{j,n,k}&\equiv\abs{\la_{j}(T_{n}(f))-\la_{j,k}^{\SL}(T_{n}(f))}, &
\eps^{\textrm{NA}}_{j,n,k}&\equiv\abs{\la_{j}(T_{n}(f))-\la_{j,k}^{\textrm{NA}}(T_{n}(f))},\\
\eps^{\textrm{MNA}}_{j,n,k}&\equiv\abs{\la_{j}(T_{n}(f))-\la_{j,k}^{\textrm{MNA}}(T_{n}(f))}, &
\eps^{\textrm{NAS}}_{j,n,k}&\equiv\abs{\la_{j}(T_{n}(f))-\la_{j,k}^{\textrm{NAS}}(T_{n}(f))},
\end{align*}}

and the respective maximum absolute errors

{\small\begin{align*}
\eps^{\SL}_{n,k}&\equiv\max\{\eps^{\SL}_{j,n,k}\colon j=1,\ldots,n\}, &
\eps^{\textrm{NA}}_{n,k}&\equiv\max\{\eps^{\textrm{NA}}_{j,n,k}\colon j=1,\ldots,n\}],\\
\eps^{\textrm{MNA}}_{n,k}&\equiv\max\{\eps^{\textrm{MNA}}_{j,n,k}\colon j=1,\ldots,n\}, &
\eps^{\textrm{NAS}}_{n,k}&\equiv\max\{\eps^{\textrm{NAS}}_{j,n,k}\colon j=1,\ldots,n\}.
\end{align*}}

We start with an example involving a well-known simple-loop symbol for which we can exactly calculate the coefficients $r_{k}$ in \eqref{eq:MainExp} easily, thus we will be able to compare the accuracy of the different eigenvalue approximations.

\begin{example}[A simple-loop symbol]
Consider the even simple-loop symbol given by

\begin{equation}\label{eq:KMS}
f(\tht)\coloneqq \frac{(1+\rho)^{2}}{2}\cdot\frac{1-\cos(\tht)}{1-2\rho\cos(\tht)+\rho^{2}}\qquad(0\le\tht\le2\pi),
\end{equation}
for a constant $0<\rho<1$, see Figure \ref{fg:KMS}. This symbol was inspired in the Kac--Murdock--Szeg\H{o} Toeplitz matrices introduced in \cite{KaMu53} and subsequently studied in \cite{Tr99,Tr10}, which usually are present in important physics models.

\begin{figure}[ht]
\centering
\includegraphics[width=0.4\textwidth]{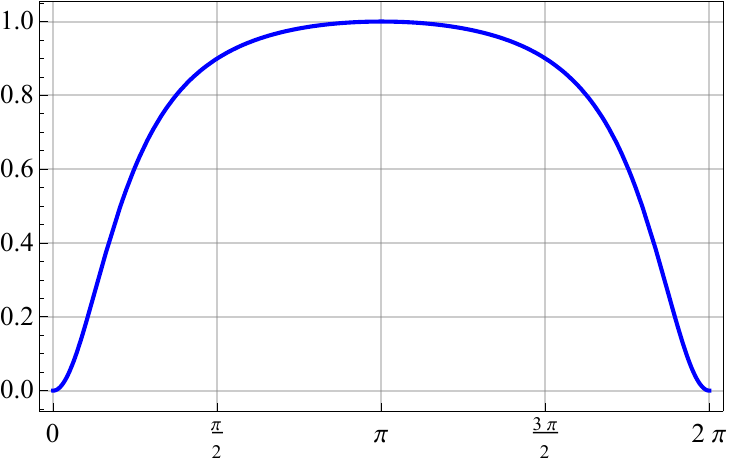}
\put(-67,90){\footnotesize $f$}
\hspace{2mm}
\includegraphics[width=0.4\textwidth]{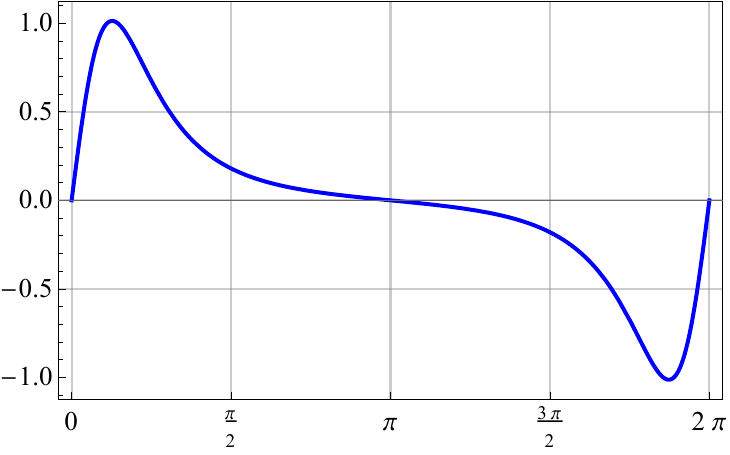}
\put(-67,90){\footnotesize $f'$}\\ \bigskip
\includegraphics[width=0.4\textwidth]{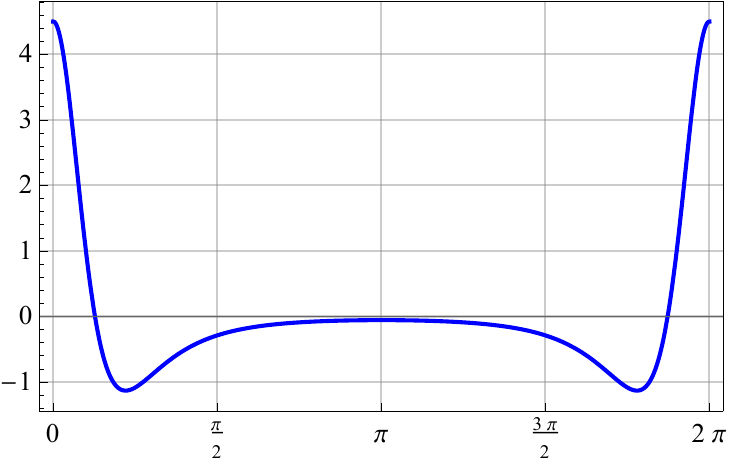}
\put(-67,90){\footnotesize $f''$}
\caption{The symbol $f$ in \eqref{eq:KMS} and its first two derivatives for $\rho=\frac{1}{2}$.}\label{fg:KMS}
\end{figure}

The respective Fourier coefficients can be calculated as $\mathfrak{a}_{k}(f)=\frac{1}{4}(\rho^{2}-1)\rho^{\,\abs{k}-1}$ for $k\ne0$ and $\frac{1}{2}(1+\rho)$ for $k=0$.
We then have
\[\norm{f}_{\al}=\frac{1+\rho}{2}+\frac{\rho^{2}-1}{2\rho}\sum_{k=1}^{\nf}\rho^{k}(k+1)^{\al},\]
which is finite for every $\al>0$, the remaining simple-loop conditions are easily verified in Figure~\ref{fg:KMS}. Then $f\in\SL^{\al}$ for any $\al>0$.
According to \cite[Sc.4]{BoSe21}, the function $\eta$ in \eqref{eq:Magic} is nicely given by
\[\eta(s)=2\arctan\Big(\frac{\rho\sin(s)}{1-\rho\cos(s)}\Big).\]

The Figure \ref{fg:ErrKMS} and the Table \ref{tb:ErrKMS} shows that, in this case, the approximation $\la_{j,k}^{\textrm{MNA}}(T_{n}(f))$ can produce good results until the level $k=2$ but it becomes unstable from this point on. While our algorithm \eqref{eq:NAS} is still producing fine results in the $4$th level and, for a matrix of size $4096$, it reaches machine-precision from level $3$.
The Figure \ref{fg:ErrKMS} reveals also, that for the symbol \eqref{eq:KMS}, our proposed algorithm \eqref{eq:NAS} can match the exact asymptotic simple-loop expansion until level $4$.
\end{example}

\begin{figure}[ht]
\centering
\includegraphics[width=0.9\textwidth]{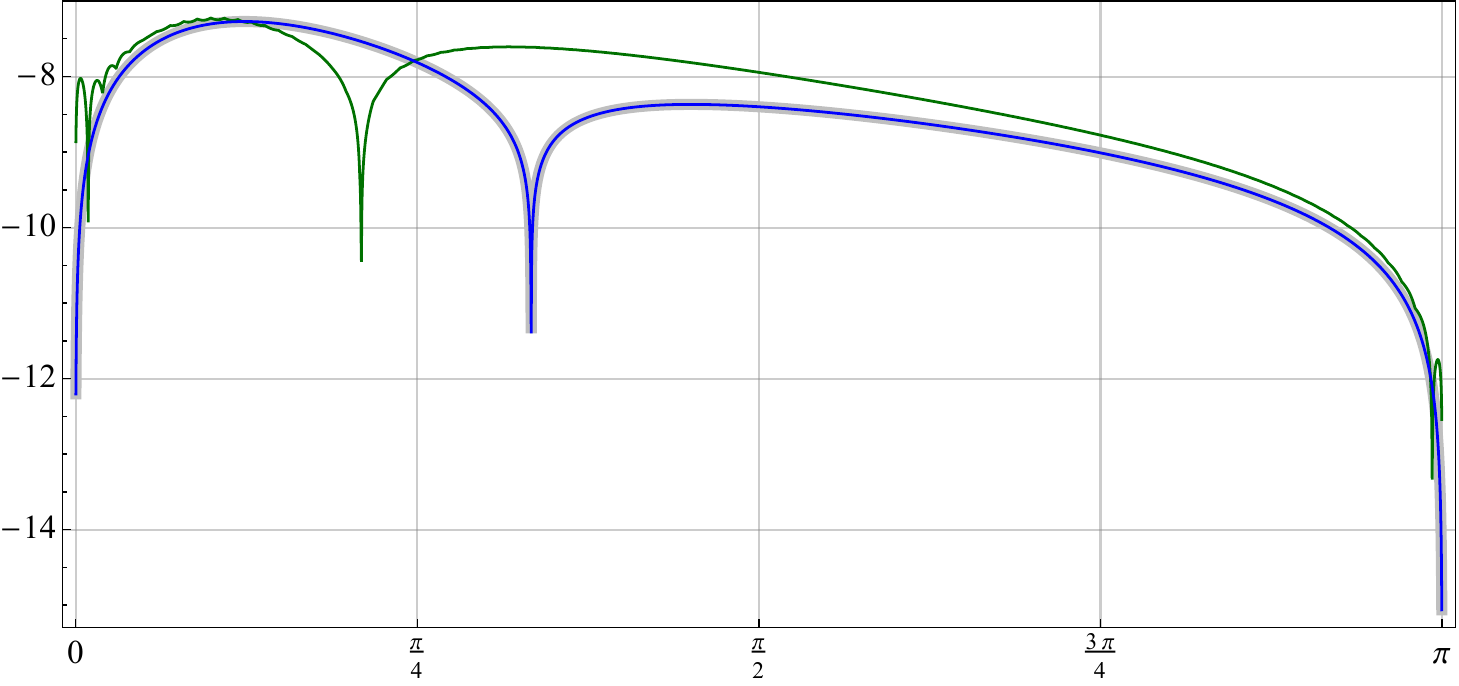}
\put(-156,145){\footnotesize $k=2$}\medskip\\
\includegraphics[width=0.9\textwidth]{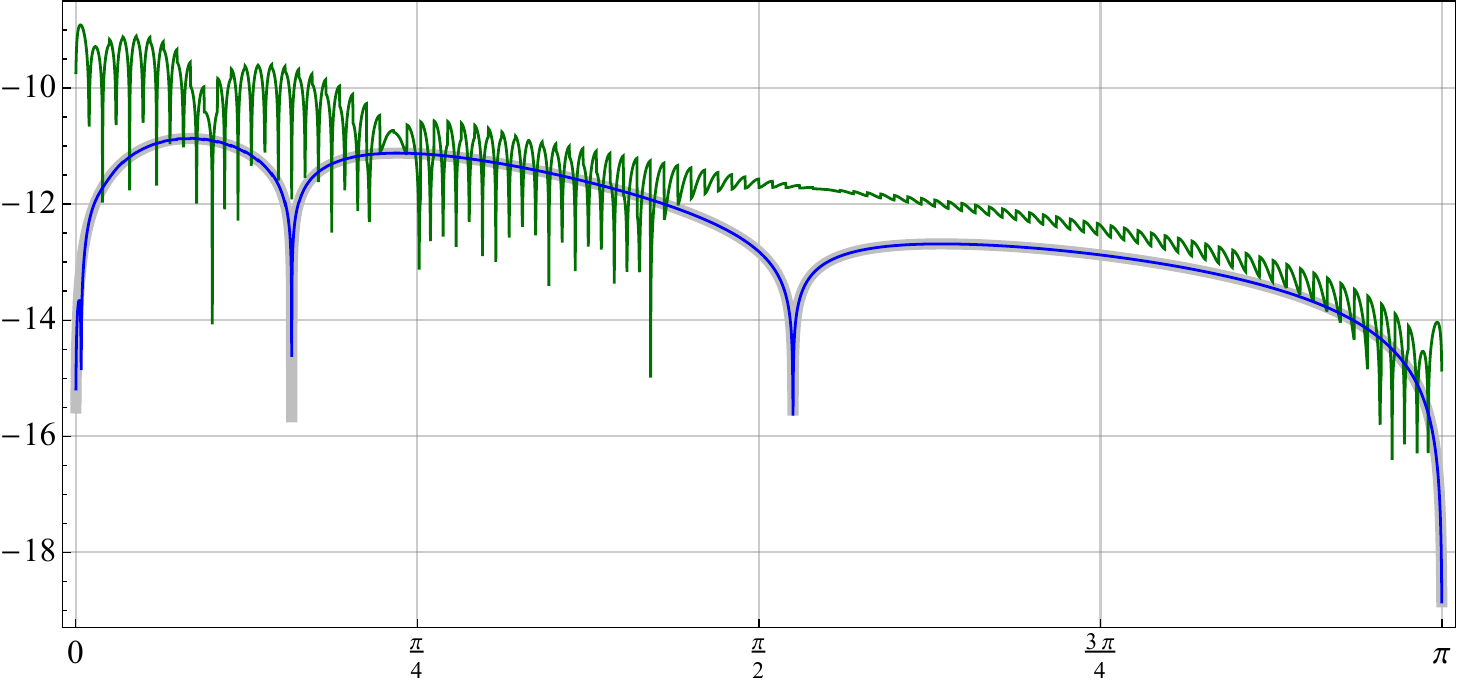}
\put(-156,145){\footnotesize $k=3$}\medskip\\
\includegraphics[width=0.9\textwidth]{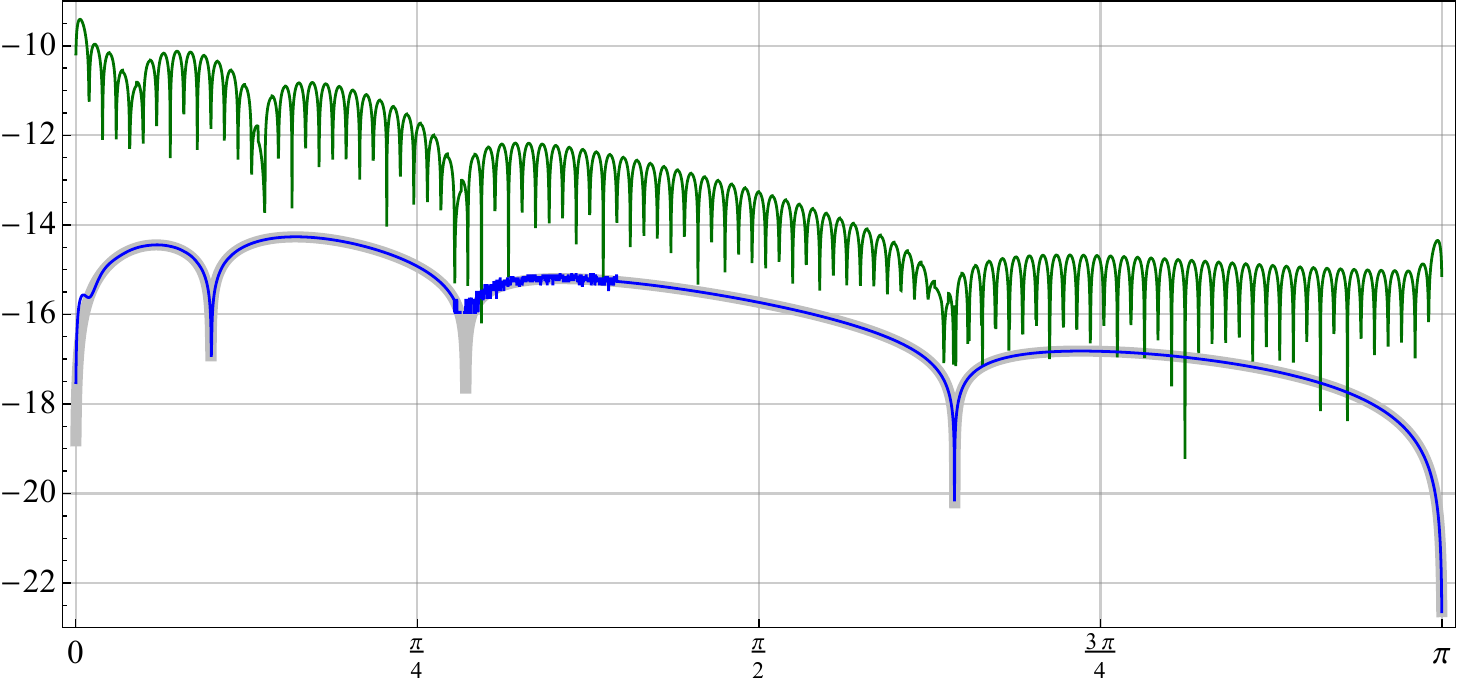}
\put(-156,145){\footnotesize $k=4$}
\caption{The base-10 logarithm for the individual errors $\eps_{j,n,k}^{\SL}$ (thick gray), $\eps_{j,n,k}^{\textrm{MNA}}$ (green), and $\eps_{j,n,k}^{\textrm{NAS}}$ (blue) for the symbol $f$ in \eqref{eq:KMS} with $\rho=\frac{1}{2}$, a matrix size $n=4096$, a grid size $n_{1}=100$, and different levels $k$. The proposed numerical algorithm \eqref{eq:NAS} is represented by the blue curve, while the exact asymptotic simple-loop expansion, by the gray curve.}\label{fg:ErrKMS}
\end{figure}
\clearpage

{\renewcommand{\arraystretch}{1.6}
\begin{table}[ht]
\centering
\caption{The maximum errors $\eps_{n,k}^{\textrm{MNA}}$, $\eps_{n,k}^{\textrm{NAS}}$, and maximum normalized errors $(n+1)^{k}\eps_{n,k}^{\textrm{NAS}}$ for the levels $k=1,2,3,4$ and different matrix sizes $n$, corresponding to the symbol \eqref{eq:KMS} with $\rho=\frac{1}{2}$. We used a grid of size $n_{1}=100$.}\label{tb:ErrKMS}
{\footnotesize\begin{tabular}{| r| l| l| l| l| l| }
\hline
\multicolumn{1}{| c| }{$n$} & \multicolumn{1}{| c| }{$256$} & \multicolumn{1}{| c| }{$512$} & \multicolumn{1}{| c| }{$1024$} & \multicolumn{1}{| c| }{$2048$} & \multicolumn{1}{| c| }{$4096$} \\ \hline\hline
$\eps_{n,1}^{\textrm{MNA}}$ & $3.0897\x10^{-3}$ & $1.5494\x10^{-3}$ & $7.7577\x10^{-4}$ & $3.8816\x10^{-4}$ & $1.9415\x10^{-3}$ \\ \hline
$\eps_{n,1}^{\textrm{NAS}}$ & $3.0897\x10^{-3}$ & $1.5494\x10^{-3}$ & $7.7577\x10^{-4}$ & $3.8816\x10^{-4}$ & $1.9415\x10^{-4}$ \\ \hline
$(n+1)\eps_{n,1}^{\textrm{NAS}}$ & $7.9405\x10^{-1}$ & $7.9482\x10^{-1}$ & $7.9517\x10^{-1}$ & $7.9534\x10^{-1}$ & $7.9542\x10^{-1}$ \\ \hline\hline
$\eps_{n,2}^{\textrm{MNA}}$ & $1.5098\x10^{-5}$ & $3.7949\x10^{-6}$ & $9.5024\x10^{-7}$ & $2.3819\x10^{-7}$ & $5.9529\x10^{-8}$ \\ \hline
$\eps_{n,2}^{\textrm{NAS}}$ & $1.3575\x10^{-5}$ & $3.4113\x10^{-5}$ & $8.5515\x10^{-7}$ & $2.1407\x10^{-7}$ & $5.3553\x10^{-8}$ \\ \hline
$(n+1)^{2}\eps_{n,2}^{\textrm{NAS}}$ & $8.9661\x10^{-1}$ & $8.9775\x10^{-1}$ & $8.9844\x10^{-1}$ & $8.9875\x10^{-1}$ & $8.9890\x10^{-1}$ \\ \hline\hline
$\eps_{n,3}^{\textrm{MNA}}$ & $8.8167\x10^{-8}$ & $1.4780\x10^{-8}$ & $4.7324\x10^{-9}$ & $2.4238\x10^{-9}$ & $1.2270\x10^{-9}$ \\ \hline
$\eps_{n,3}^{\textrm{NAS}}$ & $5.4356\x10^{-8}$ & $6.8619\x10^{-9}$ & $8.6153\x10^{-10}$ & $1.0794\x10^{-10}$ & $1.3507\x10^{-11}$ \\ \hline
$(n+1)^{3}\eps_{n,3}^{\textrm{NAS}}$ & $9.2267\x10^{-1}$ & $9.2640\x10^{-1}$ & $9.2778\x10^{-1}$ & $9.2852\x10^{-1}$ & $9.2887\x10^{-1}$ \\ \hline\hline
$\eps_{n,4}^{\textrm{MNA}}$ & $6.6888\x10^{-9}$ & $3.1948\x10^{-9}$ & $1.5778\x10^{-9}$ & $7.8461\x10^{-10}$ & $3.8983\x10^{-10}$ \\ \hline
$\eps_{n,4}^{\textrm{NAS}}$ & $3.4700\x10^{-10}$ & $2.1887\x10^{-11}$ & $1.3740\x10^{-12}$ & $8.6077\x10^{-14}$ & $5.4131\x10^{-15}$ \\ \hline
$(n+1)^{4}\eps_{n,4}^{\textrm{NAS}}$ & $1.5138\x10^{\,0}$ & $1.5158\x10^{\,0}$ & $1.5166\x10^{\,0}$ & $1.5172\x10^{\,0}$ & $1.5252\x10^{\,0}$ \\ \hline
\end{tabular}}
\end{table}}

\begin{example}[A non-simple-loop symbol]
We now test our algorithm with a Real Cosine Trigonometric Polynomial (RCTP), see \cite[Sc.1]{EkGa18}. For $\ell\in\bZ_{+}$, consider the symbol
\begin{equation}\label{eq:fl}
f_{\ell}(\tht)\equiv(2-2\cos(\tht))^{\ell},\quad\tht\in Q.
\end{equation}
The respective Fourier coefficients can be exactly calculated as $\mathfrak{a}_{k}(f_{\ell})=(-1)^{k}{{2\ell}\choose{\ell+k}}$ for $\abs{k}\le\ell$ and $\mathfrak{a}_{k}(f_{\ell})=0$ in any other case, then the respective Toeplitz matrices $T_{n}(f_{\ell})$ are banded with a band of size $2\ell+1$. The case $\ell=2$ was carefully studied by Barrera and Grudsky in \cite{BaGr17}, where they formally deduced that
\[s_{j,n}=\frac{\pi(j+1)}{n+2}+\frac{u_{1,j}}{n+2}+\frac{u_{2,j}}{(n+2)^{2}}+O(h^{3}),\]
with some bounded and continuous coefficients $u_{1,j},u_{2,j}$.
See Theorem 2.5 there. The previous expansion is slightly different from \eqref{eq:MainExp} but we was able to show that our algorithm is producing fine results in this case.

It is clear that $f_{\ell}\in W^{\al}$ for any $\ell\in\bZ_{+}$ and any $\al>0$, but $f_{\ell}$ is simple-loop only when $\ell=1$ because $f''_{\ell}(0)=0$ for $\ell\ne1$. The Figures \ref{fg:Errf2}, \ref{fg:Errf3}, and the Tables \ref{tb:Errf2}, \ref{tb:Errf3}, show the data for the cases $\ell=2,3$. Since our method is based on the simple-loop expansion \eqref{eq:MainExp}, we expected difficulties for the very first eigenvalues, corresponding to the point $\tht=0$, nevertheless, the numerical approximations are good enough for machine precision purposes.

The Figure \ref{fg:NAvMNA} shows a comparison between the individual errors $\eps_{j,n,k}^{\textrm{NA}}$, given by the eigenvalue approximation of the numerical algorithm \cite{EkGa18}, and $\eps_{j,n,k}^{\textrm{MNA}}$, which corresponds to its boundary modification given by \cite[Sc.4]{BoSe21}. Then it is clear that the modified version works better.

{\renewcommand{\arraystretch}{1.5}
\begin{table}[hb]
\centering
\caption{The maximum errors $\eps_{n,k}^{\textrm{MNA}}$, $\eps_{n,k}^{\textrm{NAS}}$, and maximum normalized errors $(n+1)^{k}\eps_{n,k}^{\textrm{NAS}}$ for the levels $k=1,2,3,4$ and different matrix sizes $n$, corresponding to the symbol \eqref{eq:fl} with $\ell=2$. We used a grid of size $n_{1}=100$.}\label{tb:Errf2}
{\footnotesize\begin{tabular}{| r| l| l| l| l| l| }
\hline
\multicolumn{1}{| c| }{$n$} & \multicolumn{1}{| c| }{$256$} & \multicolumn{1}{| c| }{$512$} & \multicolumn{1}{| c| }{$1024$} & \multicolumn{1}{| c| }{$2048$} & \multicolumn{1}{| c| }{$4096$} \\ \hline\hline

$\eps_{n,1}^{\textrm{MNA}}$ & $1.6269\x10^{-2}$ & $8.1578\x10^{-3}$ & $4.0848\x10^{-3}$ & $2.0439\x10^{-3}$ & $1.0223\x10^{-3}$ \\ \hline
$\eps_{n,1}^{\textrm{NAS}}$ & $1.6269\x10^{-2}$ & $8.1578\x10^{-3}$ & $4.0848\x10^{-3}$ & $2.0439\x10^{-3}$ & $1.0223\x10^{-3}$ \\ \hline
$(n+1)\eps_{n,1}^{\textrm{NAS}}$ & $4.1811\x10^{\,0}$ & $4.1850\x10^{\,0}$ & $4.1869\x10^{\,0}$ & $4.1878\x10^{\,0}$ & $4.1883\x10^{\,0}$ \\ \hline\hline

$\eps_{n,2}^{\textrm{MNA}}$ & $2.2934\x10^{-5}$ & $5.7602\x10^{-6}$ & $1.4434\x10^{-6}$ & $3.6126\x10^{-7}$ & $9.0367\x10^{-8}$ \\ \hline
$\eps_{n,2}^{\textrm{NAS}}$ & $2.7270\x10^{-5}$ & $6.8421\x10^{-6}$ & $1.7136\x10^{-6}$ & $4.2880\x10^{-7}$ & $1.0725\x10^{-7}$ \\ \hline
$(n+1)^{2}\eps_{n,2}^{\textrm{NAS}}$ & $1.8011\x10^{\,0}$ & $1.8006\x10^{\,0}$ & $1.8004\x10^{\,0}$ & $1.8003\x10^{\,0}$ & $1.8002\x10^{\,0}$ \\ \hline\hline

$\eps_{n,3}^{\textrm{MNA}}$ & $7.5158\x10^{-8}$ & $9.4903\x10^{-9}$ & $1.2050\x10^{-9}$ & $1.5758\x10^{-10}$ & $2.5206\x10^{-11}$ \\ \hline
$\eps_{n,3}^{\textrm{NAS}}$ & $6.9024\x10^{-8}$ & $8.6696\x10^{-9}$ & $1.0863\x10^{-9}$ & $1.3595\x10^{-10}$ & $1.7004\x10^{-11}$ \\ \hline
$(n+1)^{3}\eps_{n,3}^{\textrm{NAS}}$ & $1.1717\x10^{\,0}$ & $1.1704\x10^{\,0}$ & $1.1698\x10^{\,0}$ & $1.1695\x10^{\,0}$ & $1.1694\x10^{\,0}$ \\ \hline\hline

$\eps_{n,4}^{\textrm{MNA}}$ & $3.7838\x10^{-9}$ & $2.2540\x10^{-10}$ & $1.6954\x10^{-11}$ & $5.9426\x10^{-12}$ & $3.2321\x10^{-12}$ \\ \hline
$\eps_{n,4}^{\textrm{NAS}}$ & $2.7800\x10^{-9}$ & $1.3631\x10^{-10}$ & $7.4328\x10^{-12}$ & $4.5503\x10^{-13}$ & $5.4968\x10^{-14}$ \\ \hline
$(n+1)^{4}\eps_{n,4}^{\textrm{NAS}}$ & $1.2128\x10^{1}$ & $9.4408\x10^{\,0}$ & $8.2044\x10^{\,0}$ & $8.2044\x10^{\,0}$ & $1.5487\x10^{1}$ \\ \hline
\end{tabular}}
\end{table}}
\vspace{-2mm}
\begin{figure}[ht]
\centering
\includegraphics[width=0.9\textwidth]{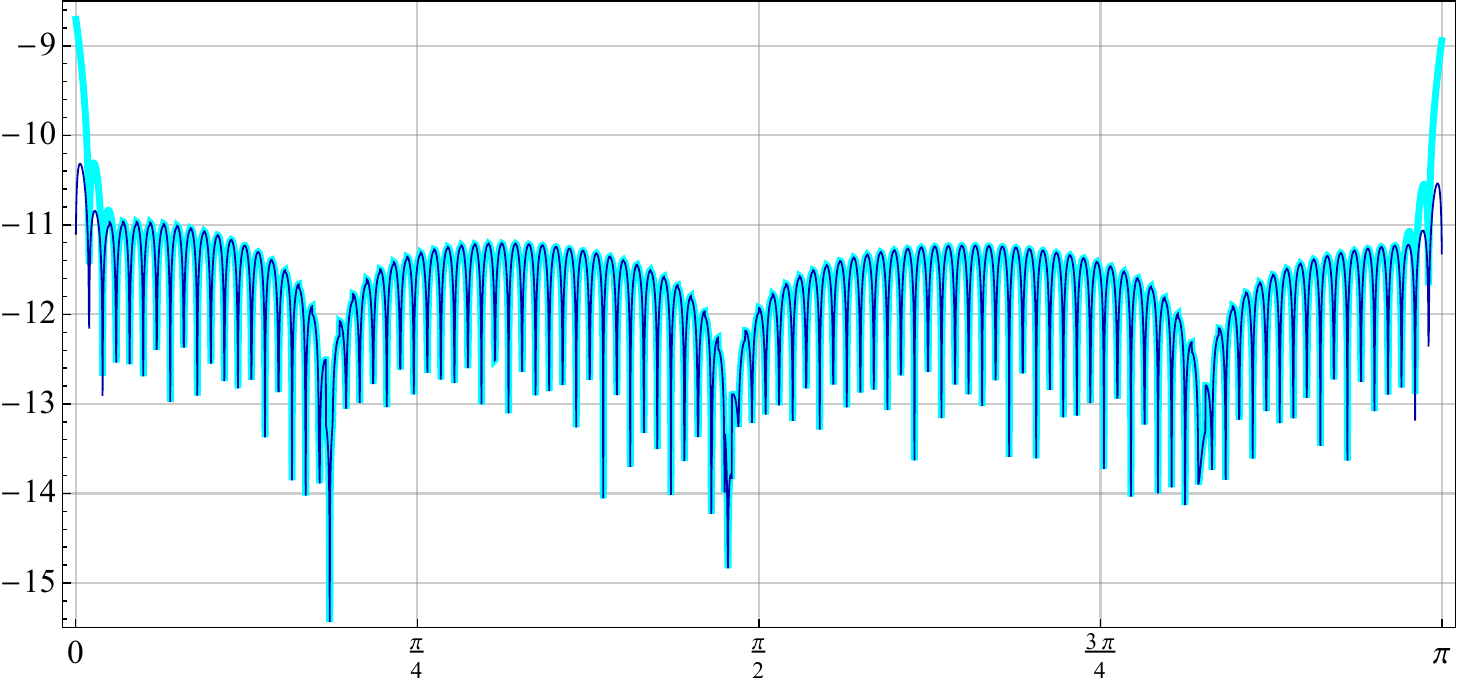}
\caption{The base-10 logarithm for the individual errors $\eps_{j,n,k}^{\textrm{NA}}$ (thick cyan curve) and $\eps_{j,n,k}^{\textrm{MNA}}$ (thin blue curve), corresponding to the algorithms \cite{EkGa18} and \cite[Sc.4]{BoSe21}, respectively. We worked with the symbol \eqref{eq:fl} where $\ell=3$, level $k=4$, a matrix size $n=4096$, and a grid size $n_{1}=100$.}\label{fg:NAvMNA}
\end{figure}

\begin{figure}[ht]
\centering
\includegraphics[width=0.9\textwidth]{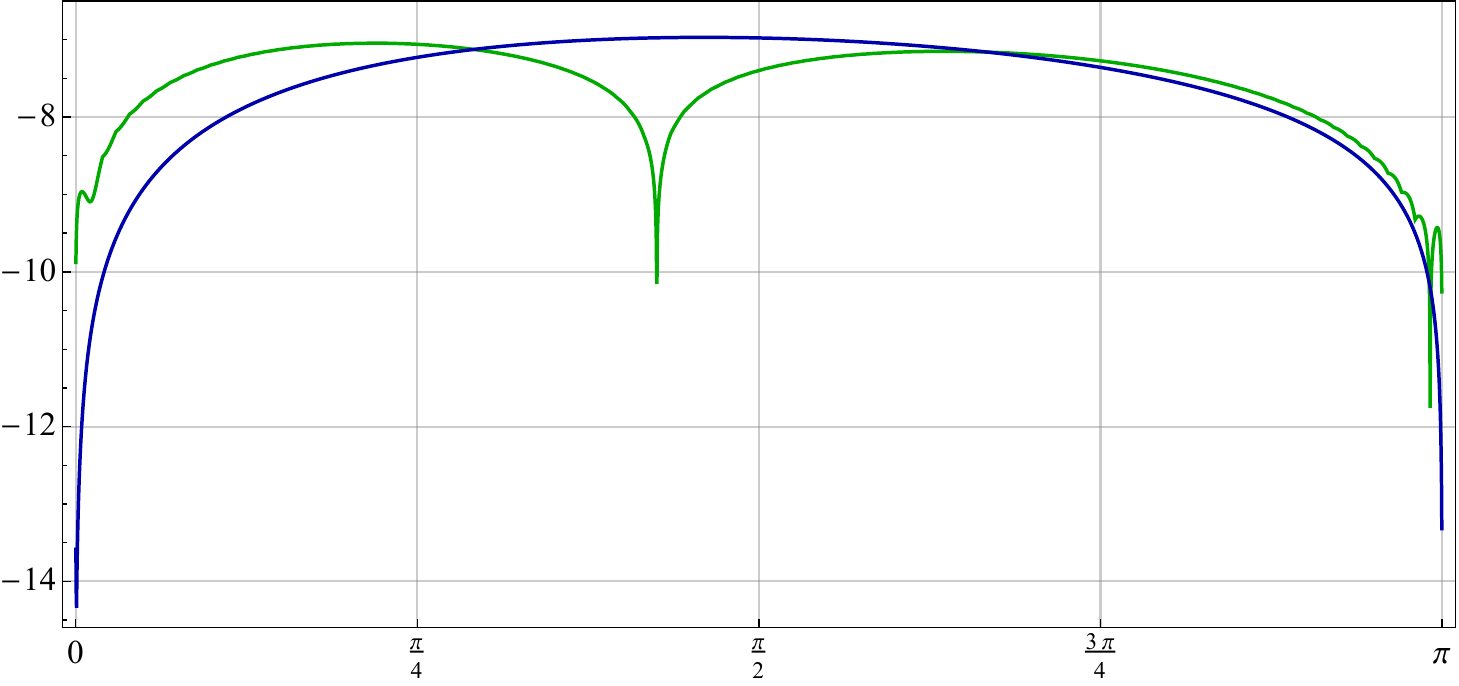}
\put(-156,145){\footnotesize $k=2$}\medskip\\
\includegraphics[width=0.9\textwidth]{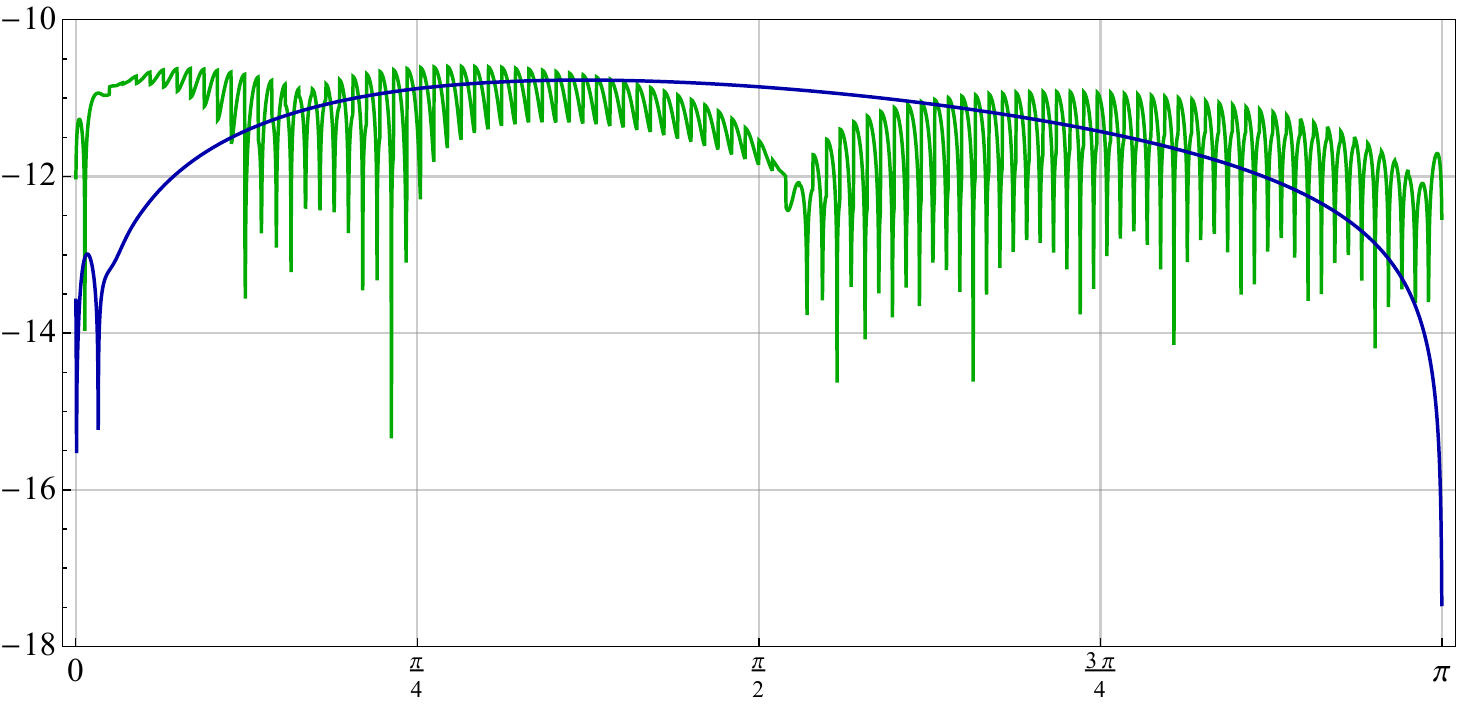}
\put(-156,145){\footnotesize $k=3$}\medskip\\
\includegraphics[width=0.9\textwidth]{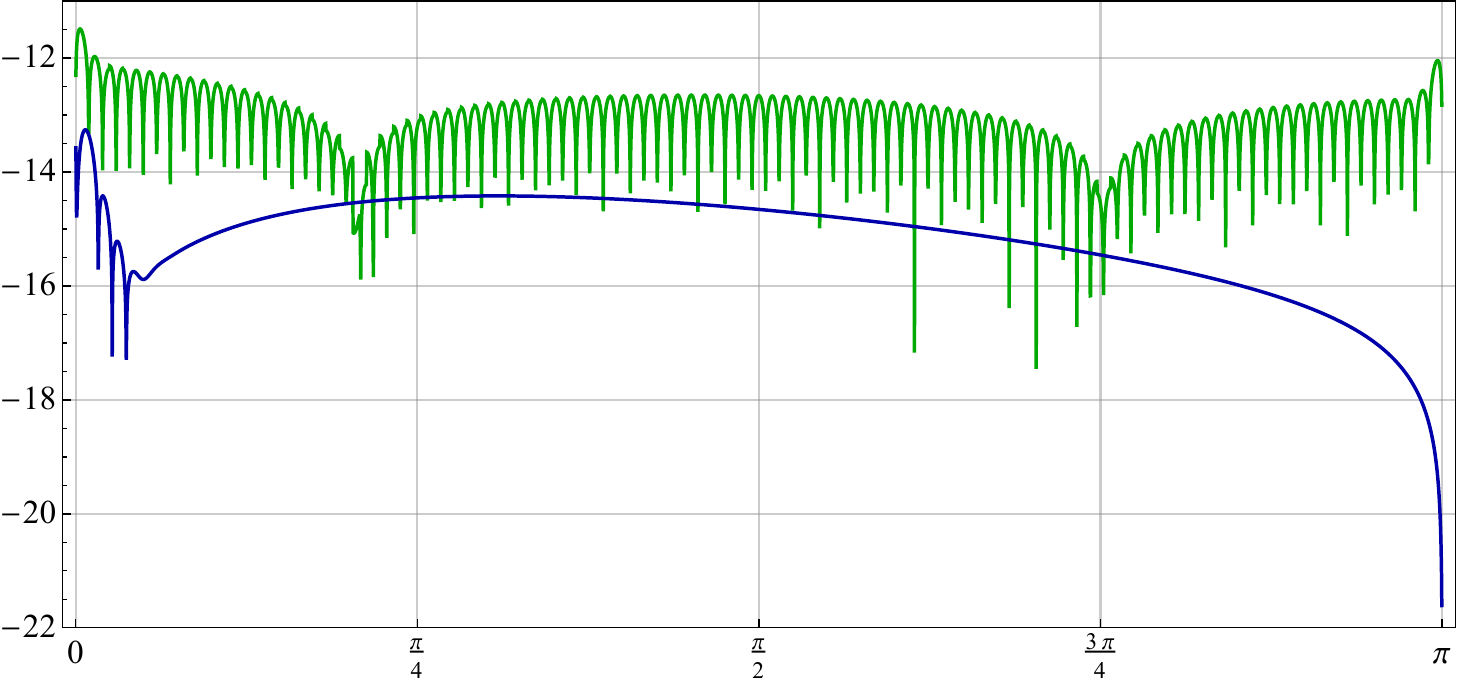}
\put(-156,145){\footnotesize $k=4$}
\caption{The base-10 logarithm for the individual errors $\eps_{j,n,k}^{\SL}$ (thick gray), $\eps_{j,n,k}^{\textrm{MNA}}$ (green), and $\eps_{j,n,k}^{\textrm{NAS}}$ (blue) for the RCTP symbol $f_{\ell}$ in \eqref{eq:fl} with $\ell=2$, a matrix size $n=4096$, a grid size $n_{1}=100$, and different levels $k$.}\label{fg:Errf2}
\end{figure}

\begin{figure}[ht]
\centering
\includegraphics[width=0.9\textwidth]{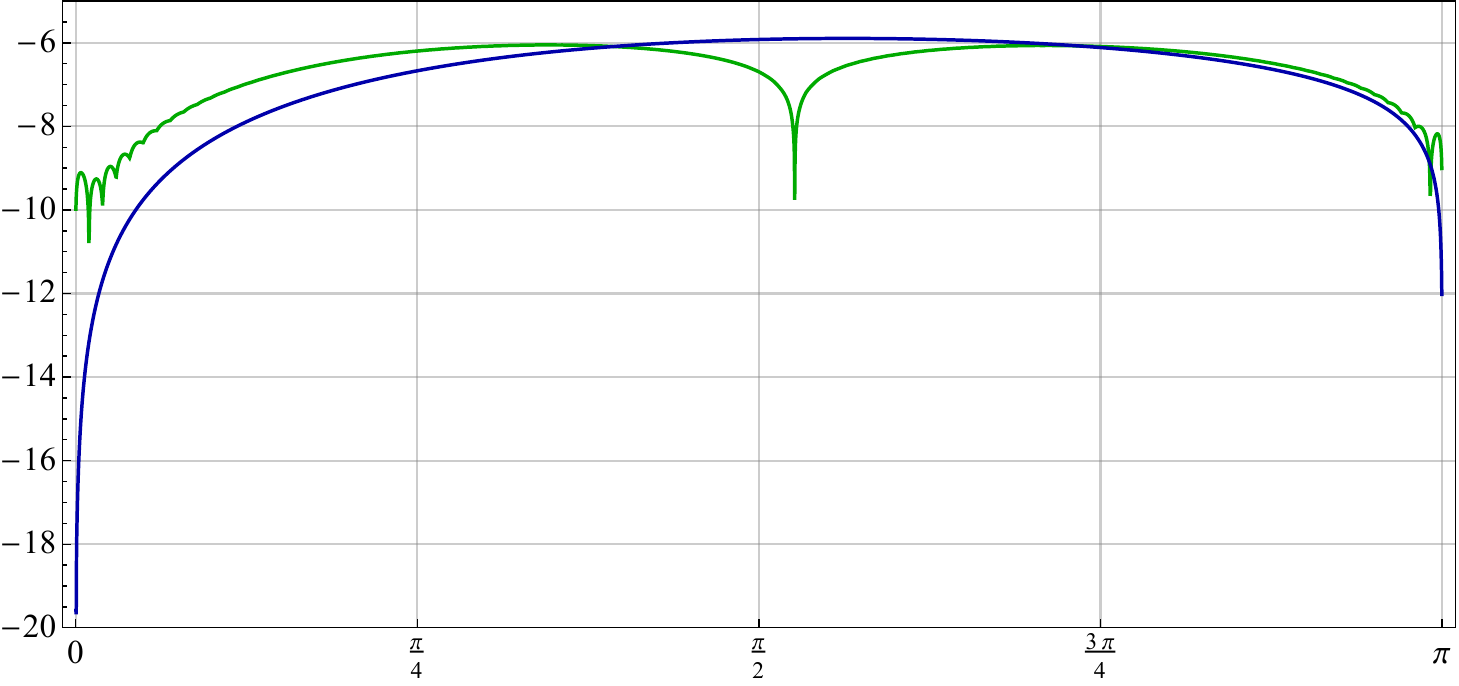}
\put(-156,145){\footnotesize $k=2$}\medskip\\
\includegraphics[width=0.9\textwidth]{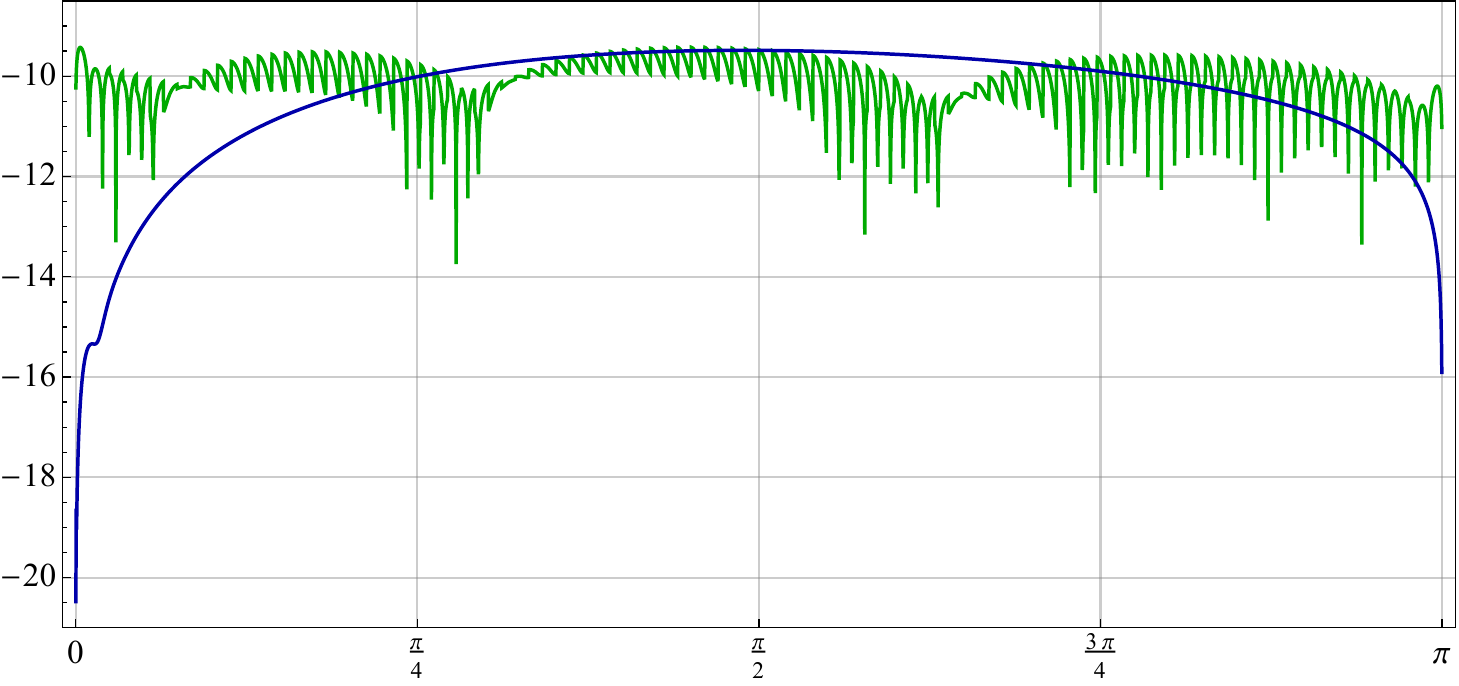}
\put(-156,145){\footnotesize $k=3$}\medskip\\
\includegraphics[width=0.9\textwidth]{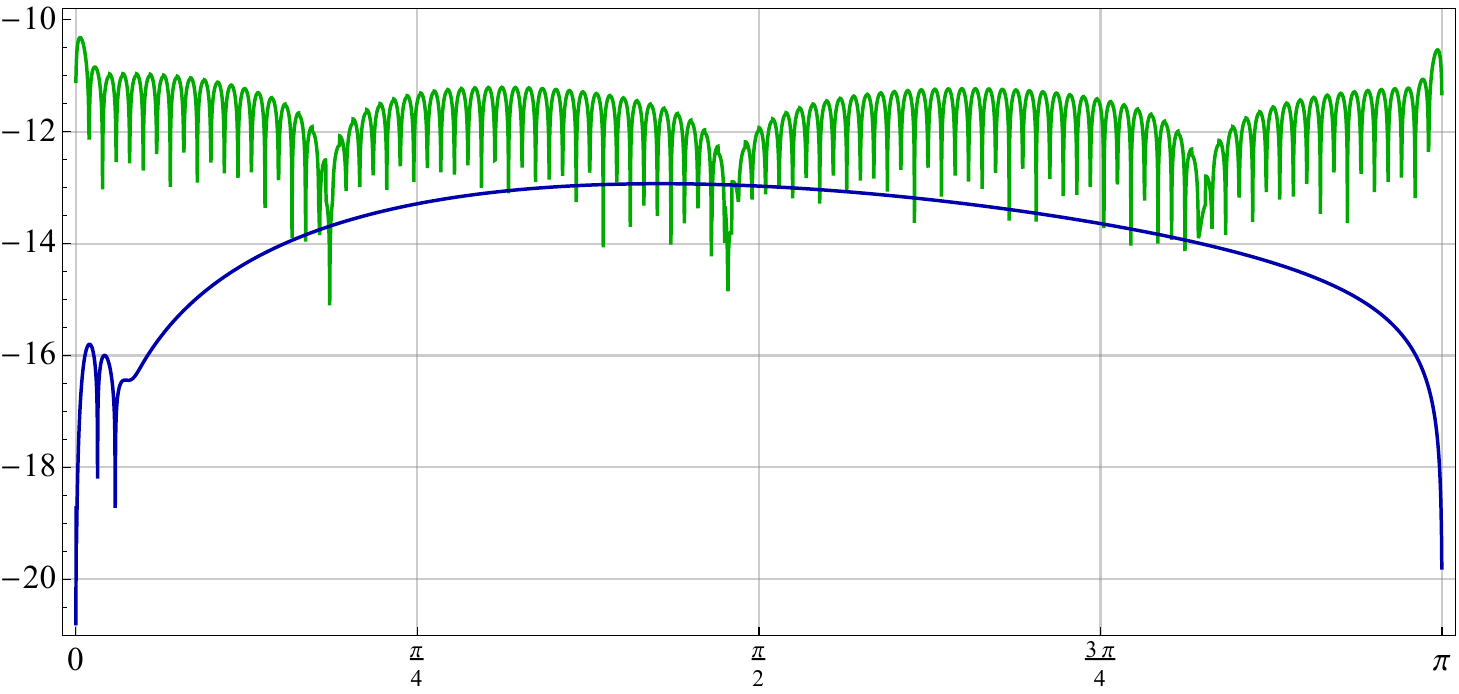}
\put(-156,145){\footnotesize $k=4$}
\caption{The same as Figure \ref{fg:Errf2} but this time with $\ell=3$.}\label{fg:Errf3}
\end{figure}
\clearpage

{\renewcommand{\arraystretch}{1.6}
\begin{table}[hb]
\centering
\caption{The same as Table \ref{tb:Errf2} but this time with $\ell=3$.}\label{tb:Errf3}
{\footnotesize\begin{tabular}{| r| l| l| l| l| l| }
\hline
\multicolumn{1}{| c| }{$n$} & \multicolumn{1}{| c| }{$256$} & \multicolumn{1}{| c| }{$512$} & \multicolumn{1}{| c| }{$1024$} & \multicolumn{1}{| c| }{$2048$} & \multicolumn{1}{| c| }{$4096$} \\ \hline\hline

$\eps_{n,1}^{\textrm{MNA}}$ & $9.1868\x10^{-2}$ & $4.6172\x10^{-2}$ & $2.3146\x10^{-2}$ & $1.1588\x10^{-2}$ & $5.7978\x10^{-3}$ \\ \hline
$\eps_{n,1}^{\textrm{NAS}}$ & $9.1868\x10^{-2}$ & $4.6172\x10^{-2}$ & $2.3146\x10^{-2}$ & $1.1588\x10^{-2}$ & $5.7978\x10^{-3}$ \\ \hline
$(n+1)\eps_{n,1}^{\textrm{NAS}}$ & $2.3610\x10^{1}$ & $2.3686\x10^{1}$ & $2.3725\x10^{1}$ & $2.3744\x10^{1}$ & $2.3753\x10^{1}$ \\ \hline\hline

$\eps_{n,2}^{\textrm{MNA}}$ & $2.2588\x10^{-4}$ & $5.6759\x10^{-5}$ & $1.4227\x10^{-5}$ & $3.5615\x10^{-6}$ & $8.9091\x10^{-7}$ \\ \hline
$\eps_{n,2}^{\textrm{NAS}}$ & $3.0497\x10^{-4}$ & $7.6550\x10^{-5}$ & $1.9176\x10^{-5}$ & $4.7989\x10^{-6}$ & $1.2003\x10^{-6}$ \\ \hline
$(n+1)^{2}\eps_{n,2}^{\textrm{NAS}}$ & $2.0143\x10^{1}$ & $2.0146\x10^{1}$ & $2.0147\x10^{1}$ & $2.0148\x10^{1}$ & $2.0148\x10^{1}$ \\ \hline\hline

$\eps_{n,3}^{\textrm{MNA}}$ & $9.5995\x10^{-7}$ & $1.2158\x10^{-7}$ & $1.5664\x10^{-8}$ & $2.1683\x10^{-9}$ & $3.7673\x10^{-10}$ \\ \hline
$\eps_{n,3}^{\textrm{NAS}}$ & $1.3355\x10^{-6}$ & $1.6765\x10^{-7}$ & $2.1002\x10^{-8}$ & $2.6281\x10^{-9}$ & $3.2868\x10^{-10}$ \\ \hline
$(n+1)^{3}\eps_{n,3}^{\textrm{NAS}}$ & $2.2669\x10^{1}$ & $2.2634\x10^{1}$ & $2.2617\x10^{1}$ & $2.2608\x10^{1}$ & $2.2604\x10^{1}$ \\ \hline\hline

$\eps_{n,4}^{\textrm{MNA}}$ & $5.7200\x10^{-9}$ & $3.9578\x10^{-10}$ & $1.7663\x10^{-10}$ & $9.3710\x10^{-11}$ & $4.8060\x10^{-11}$ \\ \hline
$\eps_{n,4}^{\textrm{NAS}}$ & $7.6467\x10^{-9}$ & $4.8020\x10^{-10}$ & $3.0083\x10^{-11}$ & $1.8824\x10^{-12}$ & $1.1772\x10^{-13}$ \\ \hline
$(n+1)^{4}\eps_{n,4}^{\textrm{NAS}}$ & $3.3358\x10^{1}$ & $3.3258\x10^{1}$ & $3.3206\x10^{1}$ & $3.3181\x10^{1}$ & $3.3168\x10^{1}$ \\ \hline
\end{tabular}}
\end{table}}
\end{example}

\begin{example}[A matrix order dependent symbol]
We now test our algorithm with a symbol which is a linear combination of RCTPs with coefficients depending on the matrix order $n$. Consider the symbol
\begin{equation}\label{eq:Fs}
F_{n}(\tht)\equiv f_{2}(\tht)+\al_{1}f_{1}(\tht)h^{2}+\al_{0}f_{0}(\tht)h^{4},
\end{equation}
where $f_{\ell}$ is given by \eqref{eq:fl} and $\al_{\ell}$ are real constants. We previously studied this symbol in \cite[Sc.4.1]{BoSe21} where we proposed an improvement to the numerical algorithm \cite{EkGa19}. This symbol commonly arises when discretizing differential equations with the Finite Differences method.

The respective Fourier coefficients can be exactly calculated using the previous example and the linearity of the Fourier transform, as $\mathfrak{a}_{k}(F_{n})=\al_{2}\mathfrak{a}_{k}(f_{2})+\al_{1}\mathfrak{a}_{k}(f_{1})h^{2}+\al_{0}\mathfrak{a}_{k}(f_{0})h^{4}$. Therefore, the Toeplitz matrices $T_{n}(F_{n})$ are banded and penta-diagonal. For implementing the numerical algorithm in \cite{EkGa19}, we need to assume an eigenvalue expansion with the form
\[\la_{j}(T_{n}(F_{n}))=f_{2}(\tht_{j,n})+\sum_{\ell=1}^{\al-1}c_{\ell}(\tht_{j,n})h^{\ell}+E_{\al,j,n},\]
where the coefficients $c_{\ell}$ are continuous and bounded functions from $[0,\pi]$ to $\bR$, and the remainder (error) term $E_{\al,j,n}$ satisfies the inequality $\abs{E_{\al,j,n}}\le\ka_{\al}h^{\al}$. For the respective boundary modification proposed in \cite{BoSe21}, we need to note that
\begin{align*}
c_{2}(0)&=\al_{1}f_{1}(0)=0, & c_{2}(\pi)&=\al_{1}f_{1}(\pi)=4\al_{1},\\
c_{4}(0)&=\al_{0}f_{0}(0)=0, & c_{4}(\pi)&=\al_{0}f_{0}(\pi)=\al_{0},
\end{align*}
while $c_{\ell}(0)=c_{\ell}(\pi)=0$ in any other case. The Figure \ref{fg:ErrF} and the Table \ref{tb:ErrF}, show the data.
\clearpage

{\renewcommand{\arraystretch}{1.6}
\begin{table}[ht]
\centering
\caption{The maximum errors $\eps_{n,k}^{\textrm{MNA}}$, $\eps_{n,k}^{\textrm{NAS}}$, and maximum normalized errors $(n+1)^{k}\eps_{n,k}^{\textrm{NAS}}$ for the levels $k=1,2,3,4$ and different matrix sizes $n$, corresponding to the symbol \eqref{eq:Fs} with $\al_{0}=3$ and $\al_{1}=2$. We used a grid of size $n_{1}=100$.}\label{tb:ErrF}
{\footnotesize\begin{tabular}{| r| l| l| l| l| l| }
\hline
\multicolumn{1}{| c| }{$n$} & \multicolumn{1}{| c| }{$256$} & \multicolumn{1}{| c| }{$512$} & \multicolumn{1}{| c| }{$1024$} & \multicolumn{1}{| c| }{$2048$} & \multicolumn{1}{| c| }{$4096$} \\ \hline\hline
$\eps_{n,1}^{\textrm{MNA}}$ & $6.0007\x10^{-3}$ & $3.0088\x10^{-3}$ & $1.5065\x10^{-3}$ & $7.5377\x10^{-4}$ & $3.7702\x10^{-4}$ \\ \hline
$\eps_{n,1}^{\textrm{NAS}}$ & $6.0007\x10^{-3}$ & $3.0088\x10^{-3}$ & $1.5065\x10^{-3}$ & $7.5377\x10^{-4}$ & $3.7702\x10^{-4}$ \\ \hline
$(n+1)\eps_{n,1}^{\textrm{NAS}}$ & $1.5422\x10^{\,0}$ & $1.5435\x10^{\,0}$ & $1.5442\x10^{\,0}$ & $1.5445\x10^{\,0}$ & $1.5446\x10^{\,0}$ \\ \hline\hline

$\eps_{n,2}^{\textrm{MNA}}$ & $3.9869\x10^{-5}$ & $1.0043\x10^{-5}$ & $2.5287\x10^{-6}$ & $6.3851\x10^{-7}$ & $1.6273\x10^{-7}$ \\ \hline
$\eps_{n,2}^{\textrm{NAS}}$ & $1.5208\x10^{-5}$ & $3.7944\x10^{-6}$ & $9.4766\x10^{-7}$ & $2.3679\x10^{-7}$ & $5.9184\x10^{-8}$ \\ \hline
$(n+1)^{2}\eps_{n,2}^{\textrm{NAS}}$ & $1.0045\x10^{\,0}$ & $9.9858\x10^{-1}$ & $9.9563\x10^{-1}$ & $9.9416\x10^{-1}$ & $9.9342\x10^{-1}$ \\ \hline\hline

$\eps_{n,3}^{\textrm{MNA}}$ & $9.2045\x10^{-8}$ & $1.0710\x10^{-8}$ & $3.6629\x10^{-9}$ & $1.8241\x10^{-9}$ & $9.1977\x10^{-10}$ \\ \hline
$\eps_{n,3}^{\textrm{NAS}}$ & $8.9731\x10^{-8}$ & $1.1313\x10^{-8}$ & $1.4203\x10^{-9}$ & $1.7792\x10^{-10}$ & $2.2264\x10^{-11}$ \\ \hline
$(n+1)^{3}\eps_{n,3}^{\textrm{NAS}}$ & $1.5231\x10^{\,0}$ & $1.5274\x10^{\,0}$ & $1.5295\x10^{\,0}$ & $1.5306\x10^{\,0}$ & $1.5311\x10^{\,0}$ \\ \hline\hline

$\eps_{n,4}^{\textrm{MNA}}$ & $3.4333\x10^{-9}$ & $2.1389\x10^{-10}$ & $5.3710\x10^{-11}$ & $2.5589\x10^{-11}$ & $1.2650\x10^{-11}$ \\ \hline
$\eps_{n,4}^{\textrm{NAS}}$ & $4.3281\x10^{-9}$ & $2.7008\x10^{-10}$ & $1.8110\x10^{-11}$ & $2.3324\x10^{-12}$ & $2.9853\x10^{-13}$ \\ \hline
$(n+1)^{4}\eps_{n,4}^{\textrm{NAS}}$ & $1.8881\x10^{1}$ & $1.8705\x10^{1}$ & $1.9990\x10^{1}$ & $4.1112\x10^{1}$ & $8.4112\x10^{1}$ \\ \hline
\end{tabular}}
\end{table}}

\begin{figure}[ht]
\centering
\includegraphics[width=0.9\textwidth]{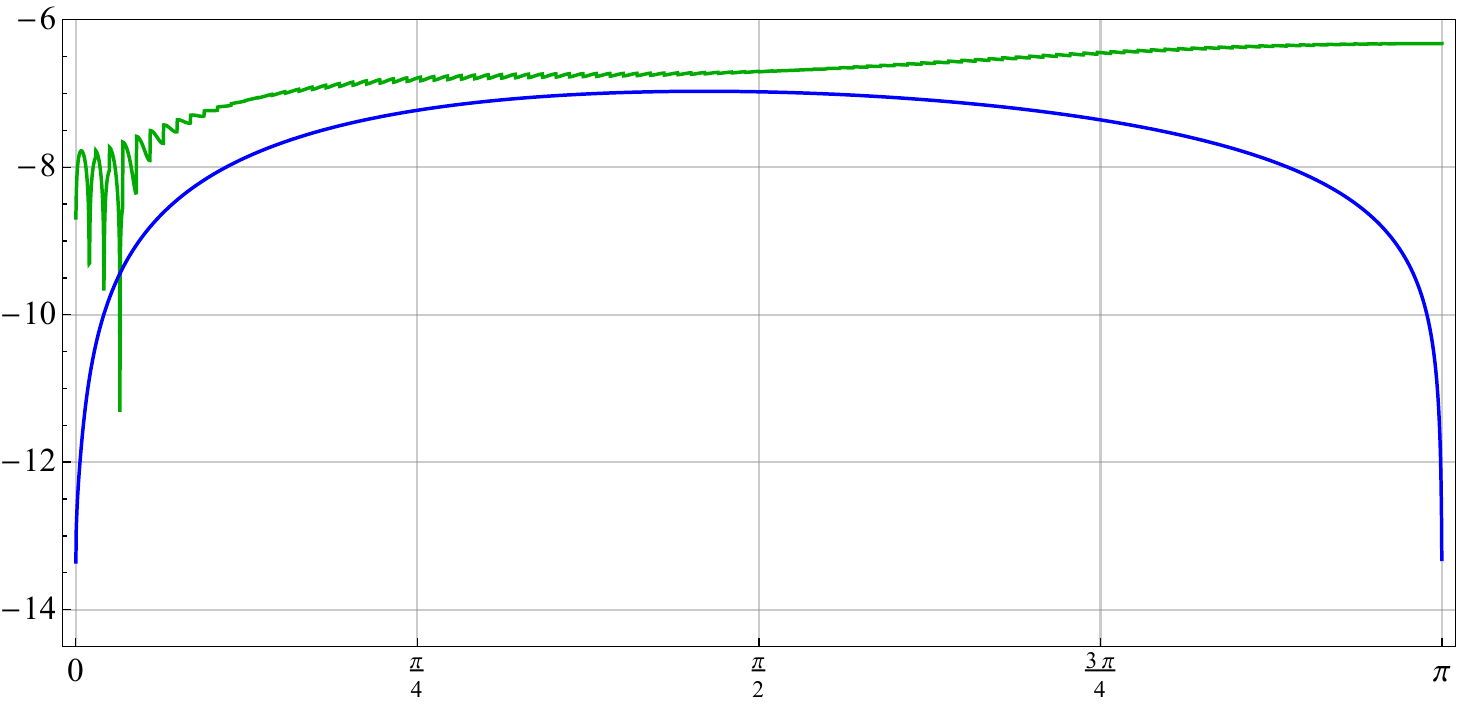}
\put(-156,145){\footnotesize $k=2$}\medskip\\
\includegraphics[width=0.9\textwidth]{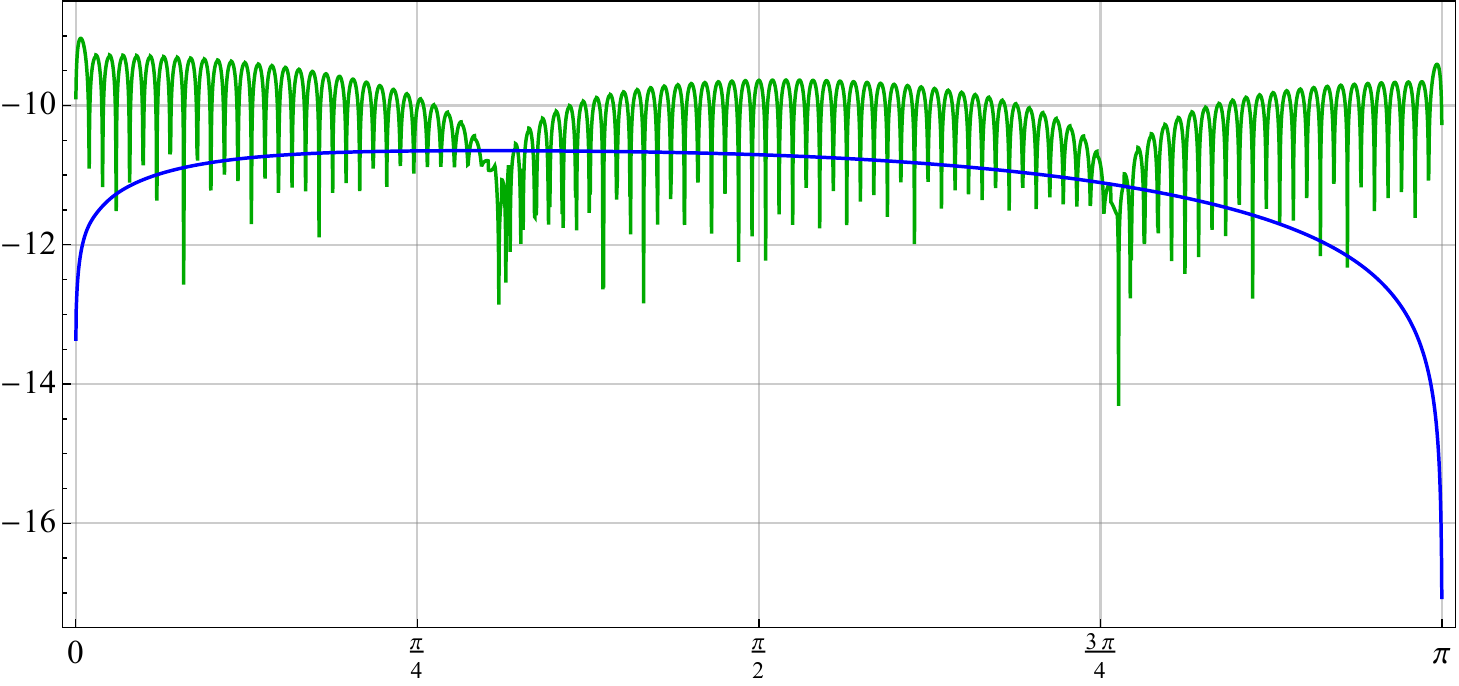}
\put(-156,145){\footnotesize $k=3$}\medskip\\
\includegraphics[width=0.9\textwidth]{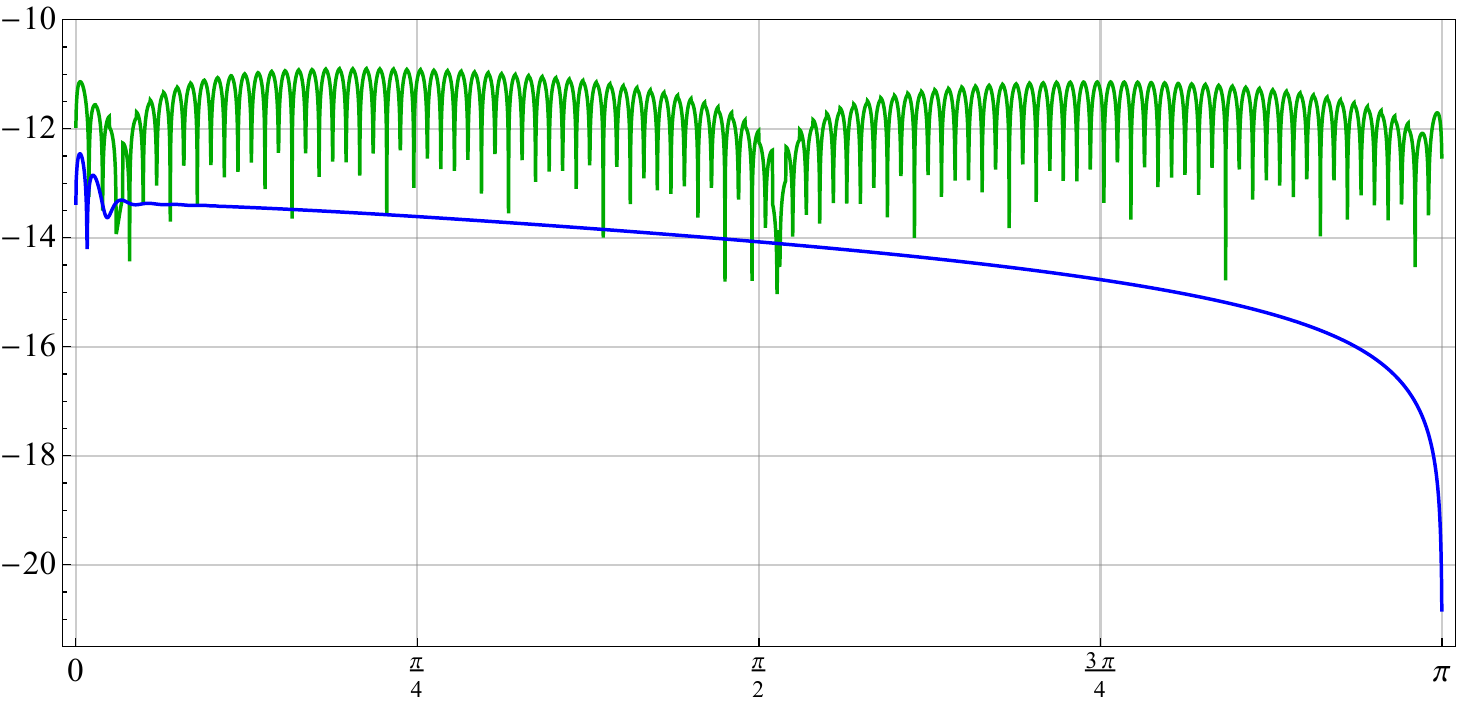}
\put(-156,145){\footnotesize $k=4$}
\caption{The base-10 logarithm for the individual errors $\eps_{j,n,k}^{\SL}$ (thick gray), $\eps_{j,n,k}^{\textrm{MNA}}$ (green), and $\eps_{j,n,k}^{\textrm{NAS}}$ (blue) for the matrix order dependent symbol $F_{n}$ in \eqref{eq:Fs}, with $\al_{0}=3$, $\al_{1}=2$, a matrix size $n=4096$, a grid size $n_{1}=100$, and different levels $k$.}\label{fg:ErrF}
\end{figure}
\end{example}

\section{Conclusions}\label{conclusions}

Under appropriate technical assumptions, the simple-loop theory allows to deduce various types of asymptotic expansions for the eigenvalues of Toeplitz matrices $T_{n}(f)$ generated by a function $f$. Independently and under the milder hypothesis that $f$ is even and monotonic over $[0,\pi]$, matrix-less algorithms have been developed for the fast eigenvalue computation of large Toeplitz matrices. These procedures work with a linear complexity in the matrix order $n$ and behind the high efficiency of such algorithms there are the expansions predicted by the simple-loop theory, combined with the extrapolation idea.

In this note we have focused our attention on a change of variable, followed by the asymptotic expansion of the new variable that is 
\[
\la_{j}(T_{n}(f))\equiv f(s_{j,n}), \qquad s_{j,n}=\tht_{j,n}+\sum_{k=1}^{\lfloor\al\rfloor}r_{k}(\tht_{j,n})h^{k}+E_{j,n,\al},
\]
and then we adapted the matrix-less to the considered new setting.
\clearpage
Numerical experiments have shown in a clear way a higher precision (till machine precision) and the same linear computation cost, when compared with the matrix-less procedures already presented in the relevant literature. More specifically, among the advantages, we concisely mention the following: 
\begin{enumerate}[a)]
\item when the coefficients of the simple-loop function are analytically known, the algorithm computes them perfectly;
\item while the proposed algorithm is better or at worst comparable to the previous ones for the computation of the inner eigenvalues, it is extremely better for the computation of the extreme eigenvalues, which are essential for determining important quantities, like the conditioning in the positive definite case.
\end{enumerate}
 
As next steps the following questions remain to be investigated:
\renewcommand\labelitemi{\tiny$\bullet$}
\begin{itemize}
\item a fine error analysis for having a theoretical explanation of the reason why the new expansion leads to a much smaller errors, when compared with the numerical results in \cite{EkFu18a,EkFu18b,EkGa18,EkGa19};
\item applications to the block cases (see \cite{BaGa20b,BaGa20a} for the theory in the block case) and related applications \cite{EkFu18a,EkFu18b} to differential problems;
\item taking inspiration from \cite{AhAl18}, extension of the technique to preconditioned case $X_{n}=T_{n}^{-1}(g)T_{n}(l)$ in which $g$ is positive over $(0,\pi)$ and is not identically constant, making use of the ergodic theorems given in \cite{Se98a}, where Theorem \ref{teoszego-tyr} is extended to the preconditioned case in a very general setting. 
\end{itemize}

\bibliography{Toeplitz}
\end{document}